\documentclass[12pt]{amsart}
\numberwithin{equation}{section}
\usepackage{amsmath,amsthm,amsfonts,amscd,eucal}

\usepackage{xypic}
\usepackage{graphicx}

\usepackage{amssymb}
\def\ca{{\mathcal A}}
\def\cb{{\mathcal B}}
\def\cc{{\mathcal C}}

\def\ce{{\mathcal E}}
\def\cf{{\mathcal F}}

\def\ch{{\mathcal H}}

\def\ck{{\mathcal K}}

\def\co{{\mathcal O}}

\def\cs{{\mathcal S}}

\def\ga{{\mathfrak A}} 
\def\gb{{\mathfrak B}}\def\gpb{{\mathfrak b}}
\def\gc{{\mathfrak C}}

\def\gf{{\mathfrak F}}
\def\gg{{\mathfrak G}}

\def\gar{{\mathfrak R}}

\def\gt{{\mathfrak T}}

\def\gz{{\mathfrak Z}}

\def\bc{{\mathbb C}}
\def\bbf{{\mathbb F}}

\def\bm{{\mathbb M}}
\def\bn{{\mathbb N}}
\def\bp{{\mathbb P}}

\def\bz{{\mathbb Z}}

\def\a{\alpha}
\def\b{\beta}
\def\g{\gamma}  \def\G{\Gamma}
\def\d{\delta}  

\def\eps{\varepsilon}

\def\l{\lambda} 

\def\m{\mu}

\def\n{\nu}
\def\r{\rho}
\def\s{\sigma} 
\def\t{\tau}
\def\f{\varphi}  \def\F{\Phi}
  \def\Th{\Theta}
\def\om{\omega} \def\Om{\Omega}

\newtheorem{thm}{Theorem}[section]
\newtheorem{lem}[thm]{Lemma}

\newtheorem{cor}[thm]{Corollary}
\newtheorem{prop}[thm]{Proposition}

\theoremstyle{definition}
\newtheorem{rem}[thm]{Remark}
\newtheorem{defin}[thm]{Definition}

\def\aut{\mathop{\rm Aut}}
\def\carf{\mathop{\rm CAR}}

\def\di{{\rm d}}

\def\ad{\mathop{\rm ad}}

\def\idd{{1}\!\!{\rm I}}
\def\alg{\mathop{\rm alg}}
\newcommand{\nn}{\nonumber}

\begin{document}

\title[exchangeable stochastic processes]
{exchangeable stochastic processes and symmetric states in quantum probability}
\author{Vitonofrio Crismale}
\address{Vitonofrio Crismale\\
Dipartimento di Matematica\\
Universit\`{a} degli studi di Bari\\
Via E. Orabona, 4, 70125 Bari, Italy}
\email{\texttt{vitonofrio.crismale@uniba.it}}
\author{Francesco Fidaleo}
\address{Francesco Fidaleo\\
Dipartimento di Matematica \\
Universit\`{a} di Roma Tor Vergata\\
Via della Ricerca Scientifica 1, Roma 00133, Italy} \email{{\tt
fidaleo@mat.uniroma2.it}}
\date{\today}

\begin{abstract}

We analyze general aspects of exchangeable quantum stochastic processes, as well as some concrete cases relevant for several applications to Quantum Physics and Probability. We establish that there is
a one--to--one correspondence between quantum stochastic processes, either preserving or not the identity, and states on free product $C^*$--algebras, unital or not unital respectively, where the exchangeable ones correspond precisely to the symmetric states.
We also connect some algebraic properties of exchangeable processes, that is the fact that they satisfy the
product state or the block--singleton conditions, to some natural ergodic ones. We then specialize the investigation for the
$q$--deformed Commutation Relations, $q\in(-1,1)$ (the case $q=0$ corresponding to
the reduced group
$C^{*}$--algebra $C^*_r(\bbf_\infty)$ of the free group on infinitely many generators), and the Boolean ones.
We also provide a generalization of De Finetti Theorem to the Fermi CAR algebra (corresponding
to the $q$--deformed Commutation Relations with $q=-1$), by showing that any state
is symmetric if and only if it is conditionally independent and identically distributed with respect to the tail algebra. The Boolean stochastic processes provide examples for which 
the condition to be independent and identically distributed w.r.t. the tail algebra, without mentioning the {\it a--priori} existence of a preserving conditional expectation, is in general meaningless in the quantum setting.
Finally, we study the ergodic properties of a class of remarkable states on the group $C^{*}$--algebra
$C^*(\bbf_\infty)$, that is the so--called Haagerup states.
\end{abstract}

\maketitle

\section{introduction}
\label{sec1}

The study of random systems with distributional symmetries, started by De Finetti in \cite{DeF}
for sequences of $2$--point valued exchangeable random variables, has known, throughout the years,
an increasing attention in many branches of Mathematics and Physics, especially in Probability Theory,
Operator Algebras, Quantum Information Theory and Entanglement. In particular, characterizing systems
of exchangeable, or symmetrically dependent, random variables is a problem of major interest, since
generally there appear strong and useful relations with independence. Some general answers
to this problem were achieved in Probability Theory. For instance, in \cite{HS} it was shown the following
generalization of De Finetti Theorem: infinite sequences of exchangeable random variables distributed
on $X=E\times E\times\cdots$, $E$ being a compact Hausdorff space, are mixtures of independent
identically distributed random variables. The case of finite sequences,
useful for the applications, was considered in \cite{DF}, whereas in \cite{DF2} the authors succeeded to characterize the so--called partially exchangeable sequences as mixtures of Markov chains, under the hypothesis of recurrence.
These results were the source for many extensions to the non commutative setting. In particular, in \cite{St2} De Finetti--Hewitt--Savage Theorem
was firstly generalized to the infinite tensor product $C^*$--algebras by showing
that the symmetric states are mixture of extremal ones, the last consisting of infinite
products of a single state.
Moreover, some general properties of exchangeable stochastic processes based on a continuous index--set were studied, see e.g. \cite{AL}. We also mention the investigation the analogous of De Finetti Theorem in a setting involving quantum symmetries, see e.g. \cite{K} and the references cited therein.

Recently, it was shown in \cite{CrF} that de Finetti Theorem still holds for the Fermi $C^*$--algebra
based on the Canonical Anticommutation Relations (CAR for short). More precisely, the convex compact set of symmetric states
states on the CAR algebra, corresponding to exchangeable stochastic processes involving Fermions,
is indeed a Choquet simplex, where the extremal (i.e. ergodic with respect to the action of all the finite
permutations of indices) ones are precisely the Araki--Moriya products (cf. \cite{AM2}) of a single, necessarily even, state. Thus, any symmetric state
is a mixture of product states by a unique barycentric measure.

As a consequence of such results, it appears now natural to address the systematic investigation of the structure
of exchangeable stochastic processes in quantum setting. Our starting point in the present paper is
to establish the perhaps expected following facts. We show that there is a one--to--one
correspondence between quantum stochastic processes based on a
$C^*$--algebra $\ga$ and states on the free product $C^*$--algebra of the same algebra $\ga$. When $\ga$
has the unity and the stochastic process under consideration has the identity, we have to consider the free product
$C^*$--algebra in the category of the unital $C^*$--algebras. We also show that the exchangeable
processes correspond precisely to the states on the free product which are invariant with respect
the action of all the permutations moving only a finite number of indices. In addition, as the
unital free product $C^*$--algebra is a quotient of the free product obtained, in a natural way,
by forgetting the identity of $\ga$, such a one--to--one correspondence passes to the quotient.
This approach, based on the universal property of free product $C^*$--algebra, can be applied to several remarkable examples. Thus, quite naturally in many cases of interest,
the investigation of
stochastic processes can be achieved directly on "concrete" $C^*$--algebras, seen as the quotient. This is the case of the
infinite tensor product and the CAR algebras, both useful for applications in Quantum Statistical Mechanics,
as well as the classical (i.e. commutative) case, covered by considering directly the free Abelian product (which corresponds to the infinite tensor product of a single Abelian $C^*$--algebra).
We also mention the cases of interest in Free Probability, that is processes on the concrete $C^*$--algebras describing the
$q$--deformed Commutation Relations, where $q=0$ corresponds to the reduced group $C^*$--algebra of the free group on infinitely many generators, or, more generally, processes on the whole free group $C^*$--algebra.

Those preliminary results clarify us the following considerations. Even if, for our purposes it appears completely
natural to study symmetric states on the free product
$C^*$--algebra, the wide generality of this structure makes almost impossible to provide too general results.
For most of the cases relevant for applications, it will be enough and potentially more useful, to
study the properties of  invariance under the natural action of the group of permutations, of the stochastic processes directly on the quotient algebra. On the other hand,
it would be nevertheless of some interest the investigation of the structure and/or the properties of
some relevant classes of states naturally arising in Free Probability. Among them, we mention the class of the so--called Haagerup states as a pivotal example, see e.g. \cite{AHO, Ha}.

In the present paper we aim to cover these topics.
More in detail, after a preliminary section containing the notations and results useful in the sequel,
in Section \ref{sec3a} we briefly describe the free product of a $C^*$--algebra in the category
of non unital and unital $C^*$--algebras and prove that the action of the group of permutations can be really
extended to both cases so that it is compatible with the passage to the quotient. Then we recall the
definition of a quantum stochastic process, and in Theorem \ref{eqstost} we establish that assigning a class of
unitarily equivalent quantum stochastic processes on a $C^*$--algebra
$\ga$ indexed by the set $J$, is equivalent to give a state on its free product $C^*$--algebra $\gb_{\ga,J}$. In this picture, the exchangeable
stochastic processes are in one--to--one correspondence with symmetric states on such algebra. Passing to the unital case, it is also
shown that exchangeable identity preserving stochastic processes on unital $C^*$--algebra $\ga$ indexed by $J$,
correspond to symmetric states on the free product $C^*$--algebra $\gc_{\ga,J}$ in the category of unital $C^*$--algebras.
Moreover, using the universal property of
$\gb_\ga$, the ergodic (i.e. extremal) properties of symmetric states on it can be exploited for
studying the structure of exchangeable stochastic processes in some "concrete" $C^*$--algebras (cf. Remark \ref{refact}),
as described below.
We end this part relative to general properties of exchangeable stochastic processes with Theorem \ref{frecazo} in Section
\ref{secerg}.
Namely, we show that some algebraic properties, such as that to satisfy the product state or block--singleton
conditions, are equivalent to some natural ergodic properties enjoyed by exchangeable stochastic processes.
In addition, this result yields some general considerations (cf. the comments at the end of Sections \ref{secerg}, and \ref{sec4}) about the boundary of symmetric states, whose structure does appear extremely complex. This circumstance changes if one
takes special $C^*$--algebras as a starting point for the investigation of properties of (exchangeable) stochastic processes. Indeed in Section \ref{sec3} we show that, for Fermion algebra, a state is symmetric if and only if the corresponding process is conditionally independent and identically distributed with respect to the tail algebra, known in Statistical Mechanics as the algebra at infinity. As in the commutative case, such a result entails, for a state on the CAR algebra, the equivalence among invariance under the action
of the group of permutations, the fact that it is a mixture of independent and identically distributed product states (cf. Theorem 5.5 of \cite{CrF}),
and the property to be conditionally independent and identically distributed with respect to the tail algebra. This appears as the first case in which de Finetti Theorem, in the form including also its conditional version, is fully extended to a non commutative $C^*$--algebra. Moreover the equivalences above are
inherited in the case of infinite tensor product algebra (see Remark \ref{inften}), whereas they fail in the general
non commutative setting, see e.g. \cite{K} for details.
Section \ref{sec4} deals with the concrete $C^*$--algebra generated by $q$--Canonical Commutation Relations for $-1<q<1$. The CAR case corresponds to $q=-1$, and the reduced group
$C^*$--algebra $C^*_r(\bbf_\infty)$ of the free group on infinitely many generators is described by
the case $q=0$. We prove that the set of the symmetric states on all these algebras (i.e. for $-1<q<1$)
including $C^*_r(\bbf_\infty)$, reduces to a singleton. We show in Section \ref{sec5}
that the same is essentially true for the case arising from the Boolean Commutation Relations. 
Contrarily to the classical situation, the Boolean case clarifies that the formulation of the condition to be independent and identically distributed w.r.t. the tail algebra without mentioning the {\it a--priori} existence of a preserving conditional expectation, is meaningless in quantum case.
Section \ref{sec4} ends with the study of the ergodic properties of a class of remarkable states in
Free Probability, that is the Haagerup states (see e.g. Corollary 3.2 of \cite{Ha} for the definition).
Such states, defined on the whole group algebra $C^*(\bbf_\infty)$ and symmetric by definition, are shown to be ergodic or, equivalently, extremal (even if it is expected that they do not fill all the extreme boundary),
but their support in the bidual algebra $C^*(\bbf_\infty)^{**}$ does not belong to the center (i.e. they do not generate any natural KMS dynamics
on the von Neumann algebra generated by the GNS representation), except for the tracial case.
These results are achieved by using the above cited Theorem \ref{frecazo}.

\section{preliminaries and notations}
\label{sec2}

Throughout the section we will present and recall some known definitions and notations useful in the sequel.

Let $J$ be an arbitrary set, and $\ga$ a $C^*$--algebra. Take a family
$\{\ga_j\}_{j\in J}\subset\ga$ of $C^*$--subalgebras. With $\alg\{\ga_j\mid j\in J\}$, we denote its $*$--algebraic span in the ambient algebra $\ga$. Let  $I_k\subset J$, $k=1,2,3$ be finite subsets.
\begin{defin}
The state $\f\in\cs(\ga)$ is said to satisfy the \emph{product state condition} (see e.g. \cite{AM2}) if
$$
\f(A_1A_2)=\f(A_1)\f(A_2)\,,
$$
whenever $A_k\in\alg\{\ga_{j_k}\mid j_k\in I_k\}$, $k=1,2$, and
$\quad I_1\cap I_2=\emptyset$.

The state $\f$ satisfies the {\it block--singleton condition} (cf. \cite{ABCL}, Definition 2.2) if
$$
\f(A_1A_2A_3)=\f(A_1A_3)\f(A_2)\,,
$$
whenever $A_k\in\alg\{\ga_{j_k}\mid j_k\in I_k\}$, $k=1,2,3$, and
$(I_1\cup I_3)\cap I_2=\emptyset$.
\end{defin}
Suppose that $\{M_j\mid j\in J\}$ are von Neumann algebras acting on the same Hilbert space $\ch$. We denote with
$\bigvee_{j\in J} M_j:=\big(\bigcup_{j\in J} M_j\big)''$
the von Neumann algebra generated by the $M_j$.

A group $G$ is said to act as a group of automorphisms of $\ga$ if there exists a representation $\a: g\in G\rightarrow \a_g\in\aut(\ga)$.
We denote by $(\ga, G)$ this circumstance.
A state $\f\in\cs(\ga)$ is called $G$--invariant if $\f=\f\, \circ \a_g$ for each $g\in G$. The subset of the $G$--invariant states is denoted by $\cs_G(\ga)$. If $\ga$ is unital, it is $*$--weakly closed and its extremal points
are called {\it ergodic} states. For $(\ga, G)$ as above, and an invariant state $\f$ on $\ga$,
$(\pi_\f,\ch_\f,U_\f,\Omega_\f)$ is the GNS covariant
quadruple canonically associated to $\f$ (see e.g. \cite{BR1, T1}).
If $(\pi_\f,\ch_\f,\Omega_\f)$ is
the GNS triple associated to $\f$, the unitary
representation $U_\f$ of $G$ on $\ch_\f$ is uniquely determined by
\begin{align*}
&\pi_\f(\a_g(A))=U_\f(g)\pi_\f(A)U_\f(g)^{-1}\,,\\
&U_\f(g)\Om_\f=\Om_\f\,,\quad A\in\ga\,, g\in G\,.
\end{align*}
If $\gz_\f:=\pi_\f(\ga)''\bigwedge\pi_\f(\ga)'$ is the center of $\pi_\f(\ga)''$, $\gb_G(\f)$ denotes its fixed point
algebra under the adjoint action $\ad(U_\f)$ of $G$, i.e.
$$
\gb_G(\f):=\gz_\f\bigwedge \{U_\f(G)\}'
$$
In addition, let $s(\f)$ be the support of $\f$ in the bidual $\ga^{**}$. Then $s(\f)\in Z(\ga^{**})$ if and only
if $\Om_\f$ is separating for $\pi_\f(\ga)''$, $Z(\gb)$ being the center of any algebra $\gb$ (see \cite{SZ}, Section 10.17).
The invariant state $\f\in\cs_G(\ga)$ is said to be $G$--{\it Abelian} if all the operators
$E_\f\pi_\f(\ga)E_\f$ mutually commute. The $C^*$--dynamical system $(\ga,G)$ is $G$--Abelian if $\f$ is
$G$--abelian for each $\f\in\cs_G(\ga)$.

The group of permutations of $J$,
$\bp_J:=\bigcup\{\bp_I|I\subseteq J\, \text{finite}\}$
 is given by the permutations leaving fixed all the elements in $J$ but a finite number of them. If $J$ is countable, we sometimes denote $\bp_J$ simply as $\bp_\infty$. If $\bp_J$ acts as a group of automorphisms on the $C^*$--algebra $\ga$, $\f\in\cs(\ga)$ is called {\it symmetric} if it is
$\bp_J$--invariant. Following the notation introduced
above, in the unital case $\cs_{\bp_J}(\ga)$ and $\ce\left(\cs_{\bp_J}(\ga)\right)$ denote respectively the
convex closed subset of all the symmetric states of $\ga$, and the ergodic ones.
Let $M$ be the Cesaro Mean w.r.t. $\bp_J$,
given for a generic object $f(g)$ by
\begin{equation*}
M\{f(g)\}:=\lim_{I\uparrow J}\frac1{|\bp_I|}\sum_{g\in \bp_I}f(g)\,,
\end{equation*}
provided the l.h.s. exists in the appropriate sense, and $I\subset J$ runs over all the finite parts of $J$.
The state $\f\in \cs_{\bp_J}(\ga)$ is called \emph{weakly clustering} if
$$
M\{\f(\a_g(A)B)\}=\f(A)\f(B)\,,\quad A,B\in\ga,\,g\in\bp_J\,.
$$
In unital case, any weakly clustering state is ergodic and the converse holds true if $\f$ is $\bp_J$--Abelian, see e.g. Proposition 3.1.12 of \cite{S}.

A {\it conditional expectation} $E: \ga\rightarrow \gb$ between $C^*$--algebras $\ga$, $\gb$ is a norm--one linear projection of
$\ga$ onto $\gb$. If in addition it preserves the state $\f\in\cs(\ga)$: $\f\circ E=\f$, then $E$ is called a $\f$--{\it preserving} conditional expectation. The reader is referred to Chapter II of  \cite{Sr} or Section IX of \cite{T1},
and the references cited therein, for the main properties and the
existence conditions of conditional expectations.

Finally, we report the following Lemma useful in the sequel. Its immediate proof is omitted.
\begin{lem}
\label{permtra}
Consider a finite interval $J=[k,l]\subset\bz$. Then there exists a cycle
$\g\in\bp_\bz$ such that $[k+1,l+1]=\g(J)$.
\end{lem}

\section{exchangeable stochastic processes}
\label{sec3a}

In order to define exchangeable stochastic processes, we preliminary and briefly describe free products
of a single $C^{*}$--algebra $\ga$ in the categories of $C^{*}$--algebras and unital $C^{*}$--algebras,
provided that $\ga$ has a unit $\idd$ for the latter.
Indeed, if $J$ is an index set, the algebraic free product $\gb^{(0)}_{\ga,J}\equiv\gb^{(0)}$ in the category of the ${*}$--algebras
is given, as a vector space, by
$$
\gb^{(0)}:=\bigoplus_{n\geq1}\left(\bigoplus_{i_1\neq i_2\neq\cdots\neq i_n}
V_{i_1}\otimes V_{i_2}\otimes\cdots\otimes V_{i_n}\right)\,,
$$
where
$V_{i}=\ga$, $i\in J$. Notice that the indices
$i_1, i_2,\dots,i_n$ will appear possibly more than once in the r.h.s.\,.
The adjoint and the product in $\gb^{(0)}$ are defined in the usual way.
In fact, let $v=A_{i_1}\otimes A_{i_2}\cdots\otimes A_{i_m}$,
$w=B_{j_1}\otimes B_{j_2}\cdots\otimes B_{j_n}$ be two reduced words of length $m$ and $n$, respectively.\footnote{As a word in the generators of an algebra or a group, we mean always a reduced word without any further specification.} Then
$v^*:=A^*_{i_m}\otimes A^*_{i_2}\cdots\otimes A^*_{i_1}$, and
$$
vw:=\left\{
\begin{array}{ll}A_{i_1}\otimes A_{i_2}\cdots\otimes A_{i_m}
B_{j_1}\otimes B_{j_2}\cdots\otimes B_{j_n}\,,\,\qquad i_m=j_1\,,\\
A_{i_1}\otimes A_{i_2}\cdots\otimes A_{i_m}\otimes
B_{j_1}\otimes B_{j_2}\cdots\otimes B_{j_n}\,,\quad i_m\neq j_1\,.
\end{array} \right.
$$
It is straightforwardly seen that the product is well defined. Indeed, let $B\in\ga\subset\ca$ where
$\ca$ is the $C^*$--algebra obtained by $\ga$ by adding a unit $\idd$. We have, for $i_1\neq i_2\neq\cdots\neq i_n$,
\begin{align}
\label{bendef}
&\sum_lA^{(l)}_{i_1}\otimes A^{(l)}_{i_2}\cdots\otimes A^{(l)}_{i_m}B\\
=\bigg(\sum_lA^{(l)}_{i_1}&\otimes A^{(l)}_{i_2}\cdots\otimes A^{(l)}_{i_m}\bigg)
(\idd\otimes\idd\otimes\cdots\otimes B)\,,\nn
\end{align}
where $l$ runs over a finite set.
Thus, the adjoint and product defined above extend by linearity on the whole
$\gb^{(0)}$. Notice that $\gb^{(0)}$ does not have the unit even if $\ga$ has.

The algebraic free product $\gc_{\ga,J}^{(0)}\equiv\gc^{(0)}$ in the category of the unital ${*}$--algebras
is given as a vector space,
$$
\gc^{(0)}:=\bc\idd\oplus\bigoplus_{n\geq1}\left(\bigoplus_{i_1\neq i_2\neq\cdots\neq i_n}
W_{i_1}\otimes W_{i_1}\otimes\cdots\otimes W_{i_n}\right)\,,
$$
after writing
$\ga=\bc\idd\oplus W$ as a vector space (cf. Lemma 4.21 of \cite{Rud}), and $W_i=W$, $i\in J$. Also in this situation the indices
$i_1, i_2,\dots,i_n$ will appear possibly more than once in the r.h.s.\,.
The $*$--operation works as in $\gb^{(0)}$, whereas the definition of the product has to be done more carefully.
The identity (i.e. the zero degree word) is the neutral element for it. Concerning the other situation relative
to homogeneous words $v$, $w$ of length $n,m>0$ as above (with elementary tensors belonging to $W$), we get the
same result of the non unital case when $i_m\neq j_1$, i.e.
$$
vw:=A_{i_1}\otimes A_{i_2}\cdots\otimes A_{i_m}\otimes
B_{j_1}\otimes B_{j_2}\cdots\otimes B_{j_n}\,.
$$
If $i_m= j_1$ we uniquely write $
A_{i_m}B_{j_1}=\a\idd+\b C$,
where $\a,\b\in\bc$ and $C\in W$. Then
\begin{align*}
vw:=
&\a A_{i_1}\otimes A_{i_2}\cdots\otimes A_{i_{m-1}}\otimes B_{j_2}\cdots\otimes B_{j_n}\\
+&\b A_{i_1}\otimes A_{i_2}\cdots\otimes C
\otimes B_{j_2}\cdots\otimes B_{j_n}\,.
\end{align*}
In this case, \eqref{bendef} holds true without adding any other identity to $\ga$. Then the product is again well defined.

The free product $*$--algebras $\gb^{(0)}$ and $\gc^{(0)}$ are the universal algebras making commutative the following diagrams
\begin{equation}
\label{unni}
\xymatrix{ \ga \ar[r]^{i_j} \ar[d]_{\F_{j}} &
\gb^{(0)} \ar[dl]^\F \\
B}\qquad
\xymatrix{ \ga \ar[r]^{i_j} \ar[d]_{\F_{j}} &
\gc^{(0)} \ar[dl]^\F \\
C}\qquad j\in J\,.
\end{equation}
Here, for $j\in J$, $i_j$ is the canonical embedding of $\ga$ into $\gb^{(0)}$, or into $\gc^{(0)}$ in unital case,
$B$ and $C$ are arbitrary $*$--algebras with $C$ unital, and $\F$, $\F_j$ $*$--homomorphisms preserving
the corresponding identities in the unital case. This simply means that $\F$ is uniquely determined by the $\F_j$.
Thanks to the universal character of the free product algebra, the construction of $\gc^{(0)}$ does not depend on the splitting
$\ga=\bc\idd\oplus W$, up to isomorphisms.
>From now on, if $A_j\in V_{i_j}$ ($A_j\in W_{i_j}$ in the unital case),
we will frequently use the identification $\ga_j\sim i_j(\ga)$ without further mention.
It is not difficult to check that, if $\ga$ is unital, there is a natural quotient map (i.e. a $*$-epimorphism)
$$
\r^{(0)}:\gb^{(0)}\to\gc^{(0)}\,.
$$
Indeed, since
$$
\ga=V_i=\bc\idd+W_i\sim\bc\idd\oplus W_i\,,\quad i\in J\,,
$$
such a map is induced at all the levels of the tensor products by
$$
A\in\ga\longmapsto a\idd\oplus(A-a\idd)\in\bc\idd\oplus W
$$
whereas $A=a\idd+(A- a\idd)$.
For example, if $A\in\gb^{(0)}$ is a degree 1 element, then $A=a\idd+(A-a\idd)$, and
$$
\r^{(0)}(A)=a\idd\oplus(A-a\idd)\,.
$$
At degree 2 level, $A=A_1\otimes A_2\in V_{i_1}\otimes V_{i_2}$, with $A_j=a_j\idd+(A_j-a_j\idd)$, $j=1,2$. Then
\begin{align*}
&\r^{(0)}(A_1\otimes A_2)=a_1a_2\idd\oplus a_2(A_1-a_1\idd)\oplus(a_1(A_2-a_2\idd)\\
\oplus&(A_1-a_1\idd)\otimes(A_2-a_2\idd)\in\bc\idd\oplus W_{i_1} \oplus W_{i_2}
\oplus W_{i_1}\otimes W_{i_2}\,.
\end{align*}
The highest level formulas can be obtained in similar way, and $\r^{(0)}$ extends on $\gb^{(0)}$ by linearity.
By using the previous definitions of product and adjoint, it is also almost immediate to show that,
for $v,w\in\gb^{(0)}$, $\r^{(0)}(v^*)=\r^{(0)}(v)^*$ and $\r^{(0)}(vw)=\r^{(0)}(v)\r^{(0)}(w)$.
Namely, $\r^{(0)}$ is a surjective linear $*$--map preserving the algebraic structure, i.e. a
$*$--epimorphism.\footnote{The existence of $\r^{(0)}$ directly follows also from the universality of $\gb^{(0)}$:
in the l.h.s. of \eqref{unni} take $B=\gc^{(0)}$.
Since the map "forgetting the identity" $\ga\to\ga$ is one--to--one, $\r^{(0)}$ is an epimorphism.}

As $g(j_{1})\neq g(j_{2})\neq\cdots\neq g(j_{n})$
if
$j_1\neq j_2\neq\cdots\neq j_n$, $j_h\in J$, $g\in\bp_J$, the permutation group $\bp_J$ acts in a natural way on $\gb^{(0)}$,
and $\gc^{(0)}$ in the unital case,
by means of the algebraic morphisms $g\in\bp_J\mapsto\b^{(0)}_g$ and $g\in\bp_J\mapsto\g^{(0)}_g$, respectively. When $\ga$ is unital, we easily get
\begin{equation}
\label{ekv}
\r^{(0)}\circ\b^{(0)}_g=\g^{(0)}_g\circ\r^{(0)}\,,\quad g\in\bp_J\,.
\end{equation}
Let $X$ be a generic element of $\gb^{(0)}$, or $\gc^{(0)}$ in the unital
case. Define on $\gb^{(0)}$ and $\gc^{(0)}$, the extended--value seminorm
\begin{equation}
\label{semm}
\|X\|:=\sup\{\|\pi(X)\|\mid \pi\,\text{is a $*$--representation}\,\}\,.
\end{equation}
By definition, for $n=0$ in the unital case, $X$ is a multiple of the identity, whereas for $n\geq1$, it is a finite combinations of words $w_n$ of length $n$
which have the forms
$w_n=A_{1}\otimes A_{2}\cdots\otimes A_{n}$, where $A_j\in V_{i_j}$ ($W_{i_j}$ respectively) for any choice of indices
$i_1\neq i_2\neq\dots\neq i_n$.
The elementary computation $\|w_n\|\leq \prod_{j}\|A_{j}\|$ yields that the extended--valued seminorm \eqref{semm}
is effectively a seminorm. It is known that \eqref{semm} is indeed a norm
(see pag. 286 of \cite{CES} for $\gb^{(0)}$, or Proposition 2.3 in \cite{Av} for $\gc^{(0)}$ when $\ga$ is $\s$--finite),
but we do not use this fact in the sequel.
It is seen (cf. \cite{Av, CES}) that the enveloping $C^{*}$--algebras $\gb_{\ga,J}\equiv\gb$ and $\gc_{\ga,J}\equiv\gc$, of $\gb^{(0)}$ and $\gc^{(0)}$ respectively, are precisely the universal
{\it free product $C^{*}$--algebra} and {\it unital free product $C^{*}$--algebra} $\gc$ making commutative the following diagrams, analogous to \eqref{unni},
\begin{displaymath}
\xymatrix{ \ga \ar[r]^{i_j} \ar[d]_{\F_{j}} &
\gb \ar[dl]^\F \\
B}\qquad
\xymatrix{ \ga \ar[r]^{i_j} \ar[d]_{\F_{j}} &
\gc \ar[dl]^\F \\
C}\qquad j\in J\,.
\end{displaymath}
Here, $B$ is any $C^{*}$--algebra, $C$ any unital $C^{*}$--algebra, and the involved homomorphisms preserve the identities in the unital case.
\begin{prop}
The algebraic actions $\b^{(0)}_g$ and $\g^{(0)}_g$ uniquely extend to actions of the permutation group
$$
g\in\bp_J\mapsto\b_g\in\aut(\gb)\,,\quad g\in\bp_J\mapsto\g_g\in\aut(\gc)
$$
on the free $C^{*}$--algebras $\gb$ and $\gc$, respectively.

In addition, the projection map $\r^{(0)}$ described above uniquely extends to a
$C^{*}$--epimorphism $\r:\gb\to\gc$ from the free product $C^{*}$--algebra
$\gb$ onto the unital free product $C^{*}$--algebra $\gc$ fulfilling, for $\ga$ unital,
\begin{equation}
\label{ekv1}
\r\circ\b_g=\g_g\circ\r\,,\quad g\in\bp_J\,.
\end{equation}
\end{prop}
\begin{proof}
Fix $g\in\bp_J$ and define $\cb$ and $\cc$ the
pre--$C^{*}$--algebras obtained by taking quotient of $\gb^{(0)}$ and $\gc^{(0)}$ respectively,
with the ideal made of all the elements for which the seminorm \eqref{semm} vanishes. Since $\b^{(0)}_g$ and $\g^{(0)}_g$
are one--to--one norm preserving maps on the dense subsets $\cb$ and $\cc$ respectively,
they uniquely extend to automorphisms of the corresponding the enveloping $C^*$--algebras.
Moreover, by the universal property of $\gb$, $\r^{(0)}$ extends to a homomorphism $\r$ into $\gc$,
whose range contains $\cc$. By Corollary I.8.2 of \cite{T1}, it induces a $*$--isomorphism
of the quotient $C^*$--algebra onto the range of $\r$. Thus $\r$ is a $C^{*}$--epimorphism and it
satisfies \eqref{ekv1}, by taking into account the corresponding property \eqref{ekv}.
\end{proof}
Now we pass to recall the definition of a quantum stochastic process.
\begin{defin}
\label{dqspi}
A {\it stochastic process} labelled by the index set $J$ is a quadruple
$\big(\ga,\ch,\{\iota_j\}_{j\in J},\Om\big)$, where $\ga$ is a $C^{*}$--algebra, $\ch$ is an Hilbert space,
the $\iota_j$'s are $*$--homomorphisms of $\ga$ in $\cb(\ch)$, and
$\Om\in\ch$ is a unit vector, cyclic for  the von Neumann algebra
$M:=\bigvee_{j\in J}\iota_j(\ga)$ naturally acting on $\ch$.

The process is said to be {\it exchangeable} if, for each $g\in\bp_J$, $n\in\bn$, $j_1,\ldots j_n\in J$, $A_1,\ldots A_n\in\ga$
$$
\langle\iota_{j_1}(A_1)\cdots\iota_{j_n}(A_n)\Om,\Om\rangle
=\langle\iota_{g(j_{1})}(A_1)\cdots\iota_{g(j_{n})}(A_n)\Om,\Om\rangle.
$$
It is said to be {\it unital} if $\iota_{j}(\idd)=I$, $j\in J$,
provided that $\ga$ has the unit $\idd$.

Two stochastic processes $\big(\ga,\ch_i,\{\iota^{(i)}_j\}_{j\in J},\Om_i\big)$, $i=1,2$ based on the
same $C^{*}$--algebra $\ga$ and the same index set $J$, are said {\it (unitarily) equivalent} if there exists a
unitary operator $V:\ch_1\to \ch_2$ such that $V\Om_1=\Om_2$,
$V\iota^{(1)}_j(A)V^*=\iota^{(2)}_j(A)$, $a\in\ga$, $j\in J$.
\end{defin}
It is worth noticing that the first definition of stochastic processes for the quantum case
was given in \cite{AFL}, where the reader is referred for further details and comparison with classical
stochastic processes in the sense of Doob. Moreover,  in literature one can find a slightly different
definition (see, e.g. \cite{AFL}). Namely, the assignments of $\ch$ and $\Om$ are replaced by an arbitrary
unital $C^{*}$--algebra $\gg$ for which $\bigvee_{j\in J}\iota_j(\ga)$ is dense, and a state $\f$ on $\gg$ (compare with Remark \ref{refact} below). Then the equivalence condition for two processes on $\ga$ indexed by the same $J$ can be expressed in terms of
mixed moments agreement of the two states (see \cite{AFL}, Proposition 1.1).\footnote{Since many classical
stochastic processes are defined in terms of their finite--dimensional distributions, irrespective of the probability spaces,
this property is the transposition to the quantum case of the fact that two stochastic processes, defined on two different
probability spaces but having the same state space, are identified if they have the same finite--dimensional distributions.}
By using the GNS construction, it is straightforward to see the equivalence of the approaches (see. e.g. \cite{Sch}, Section 1.2),
and a quantum stochastic process defined as in Definition \ref{dqspi} is said to be in the {\it canonical form}.
Definition \ref{dqspi} is given only for discrete index but it can be easily generalized to other situations (cf. \cite{AL}). Here, we deal only with stochastic processes where
the index set $J$ is discrete, that is $J\sim\mathbb{N}$ or $\bz$ in countable situation.

Let a stochastic process be given. We introduce the linear forms $\f^{(0)}$ and $\psi^{(0)}$, and the
$*$--representations $\pi^{(0)}$ and $\s^{(0)}$ on the Hilbert space $\ch$,
respectively of $\gb^{(0)}$, and $\gc^{(0)}$ in the unital case. We firstly take $\psi^{(0)}(\idd):=1$, $\s^{(0)}(\idd):=I$
and define, on the linear generators of $\gb^{(0)}$ and $\gc^{(0)}$ respectively,
\begin{align*}
\f^{(0)}(A_{1}\otimes A_{2}\cdots\otimes A_{n})
:=&\langle\iota_{j_1}(A_1)\cdots\iota_{j_n}(A_n)\Om,\Om\rangle\,,\\
\psi^{(0)}(A_{1}\otimes A_{2}\cdots\otimes A_{n})
:=&\langle\iota_{j_1}(A_1)\cdots\iota_{j_n}(A_n)\Om,\Om\rangle\,,
\end{align*}
\begin{align*}
\pi^{(0)}(A_{1}\otimes A_{2}\cdots\otimes A_{n})
:=&\iota_{j_1}(A_1)\cdots\iota_{j_n}(A_n)\,,\\
\s^{(0)}(A_{1}\otimes A_{2}\cdots\otimes A_{n})
:=&\iota_{j_1}(A_1)\cdots\iota_{j_n}(A_n)\,,
\end{align*}
where $n\in \mathbb{N}$, $A_k\in V_{j_k}$ ($W_{j_k}$ respectively), $j_k\in J$ and $j_1\neq j_2\neq\cdots\neq j_n$.
The following relations hold, provided that the objects under consideration are
well--defined (i.e. the process is unital),
$$
\f^{(0)}(B)=\langle\pi^{(0)}(B)\Om,\Om\rangle,\, B\in\gb^{(0)},\quad
\psi^{(0)}(C)=\langle\s^{(0)}(C)\Om,\Om\rangle,\, C\in\gc^{(0)}
$$
In addition, in the unital case and/or when the process is exchangeable,
\begin{align}
\label{untl}
&\f^{(0)}=\psi^{(0)}\circ\r^{(0)}\,,\quad \pi^{(0)}=\s^{(0)}\circ\r^{(0)}\,;\\
&\f^{(0)}=\f^{(0)}\circ\b_g^{(0)}\,,\quad \psi^{(0)}=\psi^{(0)}\circ\b_g^{(0)}\,,\,\,\,g\in\bp_J\,.\nn
\end{align}
\begin{thm}
\label{strep}
Fix a stochastic process $\big(\ga,\ch,\{\iota_j\}_{j\in J},\Om\big)$. Then the linear forms and the
$*$-representations
$\f^{(0)}$, $\pi^{(0)}$ of $\gb^{(0)}$, and $\psi^{(0)}$, $\s^{(0)}$ of $\gc^{(0)}$ in the unital case,
uniquely extend to states and representations $\f$, $\pi$ of $\gb$, and $\psi$, $\s$ of $\gc$ in the unital case,
satisfying the analogues of \eqref{untl} when the process is unital and/or exchangeable.

The representations $(\ch, \pi, \Om)$, and $(\ch, \s, \Om)$ in the unital case, are the
GNS representations for the states $\f\in\cs(\gb)$ and $\psi\in\cs(\gc)$, respectively.
\end{thm}
\begin{proof}
Consider the diagrams
\begin{displaymath}
\xymatrix{ \ga \ar[r]^{i_j} \ar[d]_{\iota_{j}} &
\gb \ar[dl]^\F \\
M}\qquad
\xymatrix{ \ga \ar[r]^{i_j} \ar[d]_{\iota_{j}} &
\gc \ar[dl]^\F \\
M}\qquad j\in J\,,
\end{displaymath}
where, for each $j\in J$, the $i_j$'s are the canonical embeddings of $\ga$ in $\gb$, or in $\gc$ in the unital case.
By the universal property, the morphisms $\F$ making commutative the above diagrams are nothing but the extensions of
$\pi^{(0)}$, and $\s^{(0)}$ in the unital case, to representations
$$
\pi:\gb\to M\subset\cb(\ch)\,,\quad \s:\gc\to M\subset\cb(\ch)\,.
$$
Such representations are nondegenerate since $\Om$ is cyclic. Thanks to \eqref{untl}, it is almost immediate to verify that,
when the process is unital, $\pi$ factors through $\s$ by the epimorphism $\r$, that is $\pi=\s\circ\r$.

Concerning the functionals $\f$ and $\psi$, if $A\in\gb$ or $A\in\gc$, we get $\f(A)=\langle\pi(A)\Om,\Om\rangle$,
$\psi(A)=\langle\s(A)\Om,\Om\rangle$ in the unital case. Then they are positive and $\psi$
is automatically a state, as $\psi(\idd)=1$. Since $\Om$ is cyclic, $\s$ is nothing else than the
GNS representation of $\psi$. Moreover, for $\f$, we immediately get
$\|\f\|\leq1$.  Concerning the non unital case, fix an approximate identity $\{u_a\}\subset\gb$ which always exists, see e.g. \cite{T1}, Corollary I.7.5.
As $\pi$ is nondegenerate, we have
$\sup_a\pi(u_a)=I$. Thus,
$$
\|\f\|=\sup_a\f(u_a)=\sup_a\langle\pi(u_a)\Om,\Om\rangle
=\big\langle\big(\sup_a\pi(u_a)\big)\Om,\Om\big\rangle=1\,.
$$
Thus, also in this case, $\f$ is a state with GNS representation given by $\pi$.
It is matter of routine to check all the remaining properties.
\end{proof}
The next crucial result gives that a stochastic process on a $C^*$--algebra and a state on its free product
$C^{*}$--algebra are uniquely determined each other. In fact, we have the following
\begin{thm}
\label{eqstost}
There is a one--to--one correspondence between equivalence classes of stochastic processes
on $\ga$ and states on its free product $C^{*}$--algebra
$\gb$ which factorize through the unital $C^{*}$--algebra $\gc$, provided that $\ga$ has the unit
and the processes are unital. A class is given by exchangeable stochastic process if and only if the corresponding state is symmetric.
\end{thm}
\begin{proof}
Fix a stochastic process $\big(\ga,\ch,\{\iota_j\}_{j\in J},\Om\big)$ on the  $C^{*}$--algebra $\ga$.
From Theorem \ref{strep}, it follows that there corresponds a unique state $\f$ on $\gb$ with GNS representation $(\ch, \pi, \Om)$ and,
in the unital case, a unique state $\psi$ on $\gc$ with GNS representation $(\ch, \s, \Om)$, linked by the
required factorization property. Moreover, a process unitarily equivalent to the given one gives rise to the same state,
since it determines $*$--representations on $\ch$ which are unitarily equivalent to $\pi$ or $\s$. The proof of the last property is immediate.
\end{proof}
Following an established tradition in literature, we generally refer the terms \emph{exchangeable} and \emph{symmetric}, to stochastic processes and states, respectively. We also denote symbolically by $\gf$ the free product universal $C^*$--algebras, that is $\gf$ stands for $\gb$, or $\gc$ in the unital case.
%In addiction, $\d$ stands symbolically for the action $\b$ or $\g$ of the permutations $\bp_J$.

Very often, it is more convenient to think of a (class of)
stochastic process(es) as a state on a "concrete" $C^{*}$--algebra, rather than on a general free product $C^{*}$--algebra.
The idea can be borrowed from the construction of a stochastic process not in canonical form. Indeed, take $\ga$ a $C^{*}$--algebra,
equipped with a collection of $*$--homomorphisms
$\iota_j:\ga\to\gg$ of $\ga$ into another $C^{*}$--algebra $\gg$ equipped with an action $\a$ of the permutation group $\bp_J$.
Suppose further that the algebraic span
$\alg\{\iota_j(\ga)\mid j\in J\}$ is dense in $\gg$, and
$$
\a_g\circ\iota_j=\iota_{g(j)}\,,\quad j\in J\,,\,\,\, g\in \bp_J\,.
$$
Due to the universal property of $\gb_{\ga,J}\equiv\gb$, the $C^*$--free product of $\ga$, $\gg$ can be viewed as a quotient of
$\gb_{\ga,J}$. Then we have the following
\begin{rem}
\label{refact}
Any state on $\gg$ can be viewed as a class of unitarily equivalent stochastic processes for the algebra $\ga$ (i.e. as a state on
$\gb_{\ga,J}$)
factoring through $\gg$. The stochastic processes are exchangeable (i.e. the corresponding state on $\gb_{\ga,J}$ is symmetric)
if and only if the corresponding state on
$\gg$ is symmetric. Analogous considerations can be done for unital algebras, and for the corresponding unital processes.
\end{rem}
As examples of "concrete" $\gg$, we mention the Boolean (cf. Section \ref{sec5})
and the Canonical Anticommutation Relations (cf. Section \ref{sec3}) algebras, for the non unital and unital cases, respectively. Other remarkable examples are those arising from the infinite tensor product, where the classical case is contained, and the concrete algebra corresponding to the $q$--deformed Commutation Relations, including the free group reduced
$C^{*}$--algebra as a particular case.

\section{ergodic properties of exchangeable stochastic processes}
\label{secerg}

We start by reporting the noncommutative version of the definition of conditionally independent and identically distributed stochasticc processes, which is useful for our purposes.

Let a stochastic process $\big(\ga,\ch,\{\iota_j\}_{j\in J},\Om\big)$ be given, together with its corresponding state $\f$ on the free product $C^*$--algebra $\gf$ in unital or not unital case,
of $\ga$. Define the {\it tail algebra} of the process under consideration as
\begin{equation}
\label{ttai}
\gz^\perp_\f:=\bigcap_{\begin{subarray}{l}I\subset J,\,
I \text{finite} \end{subarray}}\bigg(\bigcup_{\begin{subarray}{l}K\bigcap I=\emptyset,
\\\,\,\,K \text{finite} \end{subarray}}\bigg(\bigvee_{k\in K}\iota_k(\ga)\bigg)\bigg)''\,.
\end{equation}
For some applications in the sequel, we provide the definition of conditionally independent and identically distributed process w.r.t.
$\gz^\perp_\f$.
\begin{defin}
\label{cocaind}
The stochastic process described by the state $\f\in\cs(\gf)$, is {\it conditionally independent and identically distributed}
w.r.t. the tail algebra if there exists a conditional expectation $E_\f:\bigvee_{j\in J}\iota_j(\ga)\to \gz^\perp_\f$ preserving the vector state $\langle\,{\bf\cdot}\,\Om_\f,\Om_\f\rangle$ such that,
\begin{itemize}
\item[(i)] $E_\f(XY)=E_\f(X)E_\f(Y)$, for each finite subsets $I,K\subset J$, $I\cap K=\emptyset$, and
$$
X\in\bigg(\bigvee_{i\in I}\iota_i(\ga)\bigg)\bigvee\gz^\perp_\f\,,
Y\in\bigg(\bigvee_{k\in K}\iota_k(\ga)\bigg)\bigvee\gz^\perp_\f\,;
$$
\item[(ii)] $E_\f(\iota_i(A))=E_\f(\iota_k(A))$ for each $i,k\in J$ and $A\in\ga$.
\end{itemize}
\end{defin}

The following results link together algebraic and ergodic properties for a symmetric state on the free product $C^*$--algebra.
Notice that the second equivalence in Theorem \ref{frecazo} below was achieved in Section 3 of \cite{ABCL},
where firstly the relations between block--singleton condition and ergodicity w.r.t. the shift were explored for symmetric faithful states. Here, after considering directly the action of the permutations, we drop faithfulness
since, if such condition holds true, or the state is asymptotically Abelian, block singleton and product state conditions are equivalent, see e.g. Proposition 5.1 in \cite{F22}.
\begin{thm}
\label{frecazo}
Consider a symmetric state $\f$ on the free product $C^*$--algebra $\gf$ . The following assertions hold true.
\begin{itemize}
\item[(i)] $\f$ satisfies the product state condition if and only if it is weakly clustering,
\item[(ii)] $\f$ is a block--singleton state if and only if $\pi_\f(\gf)''\bigwedge \{U_\f(\bp_J)\}'=\bc I$.
\end{itemize}
\end{thm}
\begin{proof}
(i) Suppose that $\f$ satisfies the product state condition. Consider two words $v,w\in\gf$ with support $I_v, I_w$ respectively.
If $I$ is a finite part of $J$, define
$A:=\{g\in\bp_I\mid I_v\cap I_{\a_g(w)}=\emptyset\}$, where hereafter $\a_g$ denotes $\b_g$, or $\gamma_g$ in the unital case. We get by applying the product state condition,
\begin{align*}
\quad\quad\bigg|\frac1{|\bp_I|}\sum_{g\in\bp_I}\f(v\a_g(w))&-\f(v)\f(w)\bigg|\\
\leq\bigg|\frac1{|\bp_I|}\sum_{g\in A}\f(v\a_g(w))- \f(v)\f(w)\bigg|
+&\bigg|\frac1{|\bp_I|}\sum_{g\in\bp_I\backslash A}\f(v\a_g(w))- \f(v)\f(w)\bigg|\\
=\bigg|\frac1{|\bp_I|}\sum_{g\in\bp_I\backslash A}\f(v\a_g(w))-&\f(v)\f(w)\bigg|
\leq2\|v\|\|w\|\frac{|A^c|}{|\bp_I|}\,.
\end{align*}
Taking the limit $I\uparrow J$, by Lemma 3.3 of \cite{CrF} one has that $\frac{|A|}{|\bp_I|}\to 1$, and
$\frac{|A^c|}{|\bp_I|}\rightarrow 0$,
where $A^c:= \bp_I\backslash A$. Thus $\f$ is weakly clustering.

Suppose now that $\f$ satisfies the weakly clustering condition. Let $v,w$ two words such that for their respective supports
$I_v\cap I_w=\emptyset$. Define $B$ the subset in $\bp_J$ leaving pointwise fixed all the elements in $I_v$. We have $|B|=(|I|-|I_v|)!$.
Since $\f$ is symmetric, one has
\begin{align*}
&\f(vw)=\frac1{(|I|-|I_v|)!}\sum_{g\in B}\f(\a_g(vw))=\frac1{(|I|-|I_v|)!}\sum_{g\in B}\f(v\a_g(w))\\
=&\frac{|I|!}{(|I|-|I_v|)!}\bigg(\frac1{|\bp_I|}\sum_{g\in\bp_I}\f(v\a_g(w))\bigg)
-\frac1{(|I|-|I_v|)!}\sum_{g\in\bp_I\backslash B}\f(v\a_g(w))\,.
\end{align*}
Taking the limit $I\uparrow J$, again Lemma 3.3 of \cite{CrF} gives $\frac{|B|}{|\bp_I|}=\frac{(|I|-|I_v|)!}{|I|!}\to1$,
$\frac{|B^c|}{|\bp_I|}\to0$, and consequently
$\frac{|B^c|}{|B|}=\frac{|B^c|}{|\bp_I|}\frac{|\bp_I|}{|B|}\to0$,
whereas the weakly clustering condition ensures that
$$
\frac1{|\bp_I|}\sum_{g\in\bp_I}\f(v\a_g(w))\to \f(v)\f(w).
$$
Then $|\f(vw)-\f(v)\f(w)|$ is infinitesimal and the product state condition follows.

(ii) We firstly observe that $\pi_\f(\gf)''\bigwedge \{U_\f(\bp_J)\}'$ is generated by cluster points
in the weak operator topology, of the Cauchy net of the Cesaro averages
$\bigg\{\frac1{|\bp_I|}\sum_{g\in\bp_I}\pi_\f(\a_g(A))\mid A\in\gf\bigg\}$, where as usual, $I$ is any finite part of $J$.
It is almost immediate to show that any of such a limit point is invariant, whereas the reverse inclusion easily follows by Kaplanski Density Theorem. Suppose now that $\f$ is a block--singleton state.  Fix three words $u,v,w$ and take $\xi=\pi_\f(u)\Om_\f$,
$\eta=\pi_\f(w)^*\Om_\f$. By counting the set
$B:=\{g\in I\mid (I_u\cup I_w)\cap I_{\a_g(v)}=\emptyset\}$, and arguing as in the previous part, we get for $I\uparrow J$,
$$
\frac1{|\bp_I|}\sum_{g\in\bp_I}\langle U_\f(g)\pi_\f(v)U_\f(g)^{-1}\xi,\eta\rangle\longrightarrow
\f(v)\langle\xi,\eta\rangle\,,
$$
as a consequence of the block--singleton condition and the fact that $\frac{|B|}{|\bp_I|}\to1$ and $\frac{|B^c|}{|\bp_I|}\to0$ as before.
This implies that the Cesaro averages of
$\{\pi_\f(\a_g(A))|A\in\gf\}$ converge to $\f(A)I$ in the weak operator topology.
As a consequence, $\pi_\f(\gf)''\bigwedge \{U_\f(\bp_J)\}'$ is trivial, since the cluster points of all
these possible averages generate the whole algebra.

Suppose now that $\pi_\f(\gf)''\bigwedge \{U_\f(\bp_J)\}'=\bc I$. Fix words $u,v,w$ such that their respective supports satisfy
$I_v\cap(I_u\cup I_w)=\emptyset$, and consider the set $B\subset \bp_I$ made by permutations leaving $I_u\cup I_w$ pointwise fixed. Since $\f$ is symmetric, we get
\begin{align*}
&\f(uvw)=\frac1{|B|}\sum_{g\in B}\f(\a_g(uvw))=\frac1{|B|}\sum_{g\in B}\f(u\a_g(v)w)\\
=&\frac{|\bp_I|}{|B|}\bigg(\frac1{|\bp_I|}\sum_{g\in\bp_I}\f(u\a_g(v)w)\bigg)
-\frac1{|B|}\sum_{g\in\bp_I\backslash B}\f(u\a_g(v)w)\,.
\end{align*}
Consider now any cluster point in the weak operator topology
$$
\lim_\b\bigg(\frac1{|\bp_{I_\b}|}\sum_{g\in\bp_{I_\b}}U_\f(g)\pi_\f(v)U_\f(g)^{-1}\bigg)=:\G\,,
$$
which exists by compactness. We have $\G\in\pi_\f(\gf)''\bigwedge \{U_\f(\bp_J)\}'$, as these cluster points generate such an algebra.
Moreover, by assumption $\G=\g I$, for such number $\g\in\bc$ depending on the chosen net. Take $\xi=\pi_\f(u)\Om_\f$,
$\eta=\pi_\f(w^{*})\Om_\f$. By usual arguments, $\frac{|B|}{|\bp_I|}\rightarrow 1$ and $\frac{|B^c|}{|B|}\rightarrow 0$, as
$I\uparrow J$. Then
\begin{align*}
&\left|\f(uvw)-\f(v)\f(uw)\right|\\=&\left|\lim_\b\bigg(\frac1{|\bp_{I_\b}|}\sum_{g\in\bp_{I_\b}}\langle U_\f(g)\pi_\f(v)U_\f(g)^{-1}\xi,\eta\rangle\bigg)-\f(v)\f(uw)\right|\\
=&\left|\g\langle\xi,\eta\rangle-\f(v)\f(uw)\right|=\left|\g\f(uw)-\f(v)\f(uw)\right|.
\end{align*}
By choosing $u,w$ as the empty words, we also find $\g=\f(v)$, and the block--singleton condition follows.
\end{proof}
We end the present section with the following considerations which distinguish the cases arising from Free Probability
from all the remaining cases relevant for the applications.

Let $\f$ be a symmetric state on a free product $C^*$--algebra $\gf$, or equivalently an exchangeable stochastic process on a $C^*$--algebra.
If it satisfies the block--singleton condition, then by (ii) of Theorem \ref{frecazo} and Theorem 4.3.20 in \cite{BR1}, it is extremal.
But a priori, one cannot say that the converse is true. A sufficient condition is that the support of $\f$ in the bidual is central,
that is $\Om_\f$ is a separating vector for $\pi_\f(\gf)''$. It is well--known that
the last condition cannot be satisfied for many cases relevant for the applications like the so--called ground states (cf. \cite{BR1}).
Suppose now $\f$ merely satisfies the product state condition. It is still extremal
by (i) of Theorem \ref{frecazo} and Proposition 3.1.10 of \cite{S}. But, again a priori, the converse is not automatically true.
Indeed, a sufficient condition is given by the asymptotic Abelianess or merely the $G$--Abelianess with $G\equiv\bp_J$, of the state under consideration, see Proposition 3.1.12 of \cite{S}.
It is expected that all such sufficient conditions, automatically true in the commutative setting, and true for some cases
relevant for the applications like the CAR or the tensor product algebras, are not generally satisfied in Free Probability.
As a consequence, the boundary of the convex of symmetric states in the free product $C^*$--algebra, or equivalently exchangeable
stochastic processes on a fixed $C^*$--algebra, contains states satisfying the product state condition, but might not generally filled by them.

\section{exchangeable processes on the CAR algebra}
\label{sec3}

The section deals with the conditional form of De Finetti--Hewitt--Savage Theorem on the Fermion algebra.
Indeed we show that a state on such algebra is symmetric if and only if it is
conditionally independent and identically distributed w.r.t the tail algebra. Since in \cite{CrF}
the authors have already characterized the symmetric states on CAR algebra as mixture of products of a single even state, as a result we finally
obtain that, being the three conditions above mutually equivalent, there is a perfect similarity with the classical setting, see e.g. \cite{Ka}, Section 1.1.

Let $J$ be an
arbitrary set. The {\it Canonical Anticommutation Relations} (CAR
for short) algebra over  $J$ is the $C^{*}$--algebra $\carf(J)$
with the identity $\idd$ generated by the set $\{a_j,
a^{\dagger}_j\mid j\in J\}$ (i.e. the Fermi annihilators and
creators respectively), and the relations
\begin{equation*}
a_{j}^{*}=a^{\dagger}_{j}\,,\,\,\{a^{\dagger}_{j},a_{k}\}=\d_{jk}\idd\,,\,\,
\{a_{j},a_{k}\}=\{a^{\dagger}_{j},a^{\dagger}_{k}\}=0\,,\,\,j,k\in J\,
\end{equation*}
where $\{A,B\}:=AB+BA$ is the anticommutator between $A$ and $B$.
The parity automorphism $\Th$ acts on the generators as
$$
\Th(a_{j})=-a_{j}\,,\,\,\Th(a^{\dagger}_{j})=-a^{\dagger}_{j}\,,\quad
j\in J\,
$$
and induces on $\carf(J)$ a $\bz_{2}$--grading, which gives $\carf(J)=\carf(J)_{+} \oplus\carf(J)_{-}$, where
\begin{align*}
&\carf(J)_{+}:=\{A\in\carf(J) \ | \ \Th(A)=A\}\,,\\
&\carf(J)_{-}:=\{A\in\carf(J) \ | \ \Th(A)=-A\}\,.
\end{align*}
Elements  in $\carf(J)_+$ and in $\carf(J)_-$ are called
{\it even} and {\it odd}, respectively. Moreover,
$\carf(J)=\overline{\carf{}_0(J)}$, where
$$
\carf{}_0(J):=\bigcup\{\carf(I)\mid\, I\subset
J\,\text{finite}\,\}
$$
is the dense subalgebra of the {\it localized elements}.
For the details relative to general properties of
$\carf(J)$ and the symmetric states on it we refer the reader to \cite{CrF} and the literature cited therein.
We only mention the fact that $\bp_J$ acts as a group of automorphisms on $\carf(J)$ (cf. \cite{CrF}, Prop. 3.2), and
each symmetric state is automatically even, that is $\f\circ \Theta=\f$ (cf. \cite{CrF}, Th. 4.1). Moreover, one can see that at least for countable $J$,
$$
\carf(J)\sim\overline{\bigotimes_{J}\bm_{2}(\bc)}^{C^*}\,.
$$
Such an isomorphism is established by a Jordan--Klein--Wigner
transformation, as shown in \cite{T1}, Exercise XIV. This transformation, contrarily to the action of the permutation group on the CAR and on tensor product algebras, does not intertwines the local structures. Thus, it cannot be directly used for our purposes.

Consider the state $\f\in\cs_{\bp_J}(\carf(J))$.
As it is even, the parity $\Th$ is unitarily implemented on $\ch_\f$ by a unitary $V_\f$ satisfying $V_\f\Om_\f=\Om_\f$
and $\ad(V_\f)(\pi_\f(A))=\pi_\f(\Th(A))$.
The {\it even} and the {\it odd} part of $\pi_\f(\carf(J))''$ are defined in a similar way as before by using $\ad(V_\f)$.

The tail algebra, known in Statistical Mechanics as the {\it algebra at infinity}, is defined as in \eqref{ttai}:
$$
\gz^\perp_\f:=\bigwedge_{\begin{subarray}{l}I\subset J,\,
I \text{finite} \end{subarray}}\left(\bigcup_{\begin{subarray}{l}K\bigcap I=\emptyset,
\\\,\,\,K \text{finite} \end{subarray}}\pi_\f(\carf(K))\right)''\,.
$$
In \cite{Rob} (see also pag. 156 of \cite{BR1}), it was proven that
$$
\gz^\perp_\f=\gz_\f\bigwedge\pi_\f(\carf(J)_+)''
$$
where as usual, $\gz_\f$ denotes the center.
It is easy to see that $\pi_\f(\carf(J)_+)''=\pi_\f(\carf(J))''_+$, so the tail algebra is automatically even:
\begin{equation}
\label{veve}
\gz^\perp_\f=\gz_\f\bigwedge\pi_\f(\carf(J))''_+\,.
\end{equation}
From now on, we take $J$
countable, say $J\equiv\mathbb{N}$ since we need the CAR algebra
is separable. We symbolically denote with an abuse of notation $\ga:=\carf(\bn)$, and with $\a$ the action of $\bp_\bn$.
We also report the following result (cf. Theorem 5.5 of \cite{CrF})
for the convenience of the reader, since it will be used in the sequel.
\begin{thm}
\label{csyx}
$\cs_{\bp_{\bn}}(\ga)$ is a Choquet simplex. Then for each
$\f\in\cs_{\bp_{\bn}}(\ga)$ there exists a unique
Radon probability measure $\m$ on $\cs_{\bp_{\bn}}(\ga)$
such that
$$
\f(A)=\int_{\ce(\cs_{\bp_{\bn}}(\ga))}\psi(A)\di\m(\psi)\,,\quad A\in\ga\,.
$$
\end{thm}
We recall that, in Probability Theory, one of the main ingredients to gain the equivalence between exchangeability
and conditionally independent and identically distributed condition for infinite sequences of random variables, is the
Hewitt--Savage Theorem \cite{HS}. It states that, for a given
exchangeable stochastic process, the symmetric and the tail
$\s$--algebras coincide. Here, before proving the main theorem, we firstly check that
$$
\gb_{\bp_\mathbb{N}}(\f)=\{\gz_\f\bigwedge \{U_\f(\bp_\mathbb{N})\}'\}_+=\gz^\perp_\f\,,
$$
that is $\gb_{\bp_\mathbb{N}}(\f)$, the fixed point algebra of the center, is always Abelian and even and coincides with the tail algebra. This result
can be seen as a generalization to the Fermi case of
the above cited commutative statement, where $\gb_{\bp_\mathbb{N}}(\f)$ and $\gz^\perp_\f$ are the counterparts of the symmetric and tail $\s$--algebras, respectively.
Furthermore, it is crucial for our purpose to show that the
map $\F_\f:\pi_\f(\ga)''\to\gb_{\bp_\mathbb{N}}(\f)$ given in Theorem 3.1 of
\cite{St1}, is precisely the conditional expectation onto the tail
algebra preserving the vector state $\langle\,{\bf\cdot}\,\Om_\f,\Om_\f\rangle$. Notice that the
existence of such a conditional expectation is not a priori guaranteed if
the corresponding state is not a trace.
\begin{lem}
\label{invper}
For each $\f\in\cs_{\bp_\mathbb{N}}(\ga)$, we get
$$
\gz^\perp_\f\subset\pi_\f(\ga)''\bigwedge U_\f(\bp_\mathbb{N})'\,.
$$
\end{lem}
\begin{proof}
For the finite set $I\subset\mathbb{N}$, define
$$
M_I:=\left(\bigcup_{\begin{subarray}{l}K\bigcap I=\emptyset,
\\\,\,\,K \text{finite} \end{subarray}}\pi_\f(\carf(K))\right)''.
$$
If $g\in\bp_I$, then $\ad(U_\f(g))(X)=X$ for each element $X\in M_I$. This
implies $\gz^\perp_\f\subset U_\f(\bp_I)'$, for $I\subset \mathbb{N}$, $I$ finite. As the l.h.s. does not depend on $I$, we get
$$
\gz^\perp_\f\subset\bigg(\bigcup_{I\subset\mathbb{N},\, I \text{finite}}U_\f(\bp_I)\bigg)'=U_\f(\bp_\mathbb{N})'\,.
$$
\end{proof}
\begin{prop}
\label{havest}
Let $\f\in\cs_{\bp_{\bn}}(\ga)$. Then
$\gz_\f^\perp=\gb_{\bp_\mathbb{N}}(\f)$, and the map $\F_\f$ is the
conditional expectation of $\pi_\f(\ga)''$ onto $\gz_\f^\perp$
preserving the vector state $\langle\,{\bf\cdot}\,\Om_\f,\Om_\f\rangle$. It assumes the form
$$
\F_\f(X)=\int^\oplus_{\ce(\cs_{\bp_{\bn}}(\ga))}\langle
X_\psi\Om_\psi,\Om_\psi\rangle\idd_{\ch_\psi}\di\m(\psi)\,, \quad
X\in\pi_\f(\ga)''\,.
$$
\end{prop}
\begin{proof}
Let $\f=\int_{{\ce(\cs_{\bp_{\bn}}(\ga))}}\psi\di\m(\psi)$ be the ergodic decomposition
of $\f\in\cs_{\bp_{\bn}}(\ga)$ given in Theorem \ref{csyx}. By Theorem 4.4.3
of \cite{BR1}, Proposition 3.1.10 of \cite{S} and Theorem IV.8.25 of \cite{T1}, we get
$$
\pi_\f(\ga)'\bigwedge\{U_\f(\bp_J)\}'=\int^\oplus_{\ce(\cs_{\bp_{\bn}}(\ga))}\bc\idd_{\ch_\psi} \di\m(\psi)
\sim L^\infty(\ce(\cs_{\bp_{\bn}}(\ga)),\m)\,.
$$
Let $x,y\in\ch_\f$ and $A\in\ga$. Then by Lebesgue Dominated
Convergence Theorem and Propositions 4.3 and 5.4 of \cite{CrF},
we obtain, for a sequence $\{I_n\}$ of finite subsets invading $\mathbb{N}$,
\begin{align*}
&\langle\F_\f(\pi_\f(A))x,y\rangle
=\lim_{I_n\uparrow\bn}
\frac1{|\bp_{I_n}|}\sum_{g\in\bp_{I_n}}\langle U_\f(g)\pi_\f(A)U_\f(g)^{-1}x,y\rangle\\
=&\lim_{I_n\uparrow\bn}
\frac1{|\bp_{I_n}|}\sum_{g\in\bp_{I_n}}\int_{\ce(\cs_{\bp_{\bn}}(\ga))}\langle U_\psi(g)\pi_\psi(A)U_\psi(g)^{-1}
x_\psi,y_\psi\rangle\di\m(\psi)\\
=&\int_{\ce(\cs_{\bp_{\bn}}(\ga))}\bigg(\lim_{I_n\uparrow\bn}
\frac1{|\bp_{I_n}|}\sum_{g\in\bp_{I_n}}\langle U_\psi(g)\pi_\psi(A)U_\psi(g)^{-1}
x_\psi,y_\psi\rangle\bigg)\di\m(\psi)\\
=&\int_{\ce(\cs_{\bp_{\bn}}(\ga))}\psi(A)
\langle x_\psi,y_\psi\rangle\di\m(\psi)
=\left\langle\left(\int^\oplus_{\ce(\cs_{\bp_{\bn}}(\ga))}\psi(A)\idd_{\ch_\psi} \di\m(\psi)\right)x,y\right\rangle\,.
\end{align*}
Then we immediately conclude that
\begin{equation}
\label{dddn}
\F_\f(\pi_\f(\ga)'')\subset\pi_\f(\ga)'\bigwedge\{U_\f(\bp_\mathbb{N})\}'\,.
\end{equation}
Consider the algebraic span $\ca:=\alg\{\pi_\f(\ga),\F(\pi_\f(\ga))\}$. If
$$
X:=\sum_k\pi_\f(A_k)\F(\pi_\f(B_k))\in\ca\,,
$$
then
\begin{equation}
\label{dddn1}
\F_\f(X)=\int^\oplus_{\ce(\cs_{\bp_{\bn}}(\ga))}f_X(\psi)\idd_{\ch_\psi} \di\m(\psi)\,,
\end{equation}
where
$$
f_X(\psi):=\sum_k\psi(A_k)\psi(B_k)\,.
$$
It is then straightforward to see that the set of the functions
$\{f_X\mid X\in\ca\}\subset C(\ce(\cs_{\bp_{\bn}}(\ga)))$ is an algebra of
continuous functions containing the constants and separating the
points of the weak--$*$ compact space $\ce(\cs_{\bp_{\bn}}(\ga))$, which is then
dense by the Stone Theorem. As a consequence, by means of
\eqref{dddn} and \eqref{dddn1} one has,
$$
\gb_{\bp_\mathbb{N}}(\f)=\F_\f(\pi_\f(\ga)'')=\pi_\f(\ga)'\bigwedge\{U_\f(\bp_\mathbb{N})\}'\subset\gz_\f\,,
$$
where the first equality and the last inclusion follow by taking into account Theorem 3.1 of \cite{St1}.
Moreover, the inclusion $\pi_\f(\ga)'\bigwedge\{U_\f(\bp_\mathbb{N})\}'\subset\gz_\f$ gives (cf. Lemma 8.4.1 of \cite{DX0}),
$$
\pi_\f(\ga)''=\int^\oplus_{\ce(\cs_{\bp_{\bn}}(\ga))}\pi_\psi(\ga)''\di\m(\psi)\,.
$$
Thus, since $\ce(\cs_{\bp_{\bn}}(\ga))$ is made of
factor states (cf. Proposition 5.7 of \cite{CrF}), the
decomposition in Theorem \ref{csyx} is indeed the factor
decomposition of $\f$ and, by Corollary IV.8.20 of \cite{T1}, $\gz_\f$ is the diagonal algebra, that is
it coincides with $\pi_\f(\ga)'\bigwedge\{U_\f(\bp_\mathbb{N})\}'$. In addition,  as $\f$ and all the
$\psi\in\ce(\cs_{\bp_{\bn}}(\ga))$ are
even, we can consider their covariant
implementations $V_\f$, and $V_\psi$ respectively, in the GNS representation. The
unitary $V_\f$ can be simultaneously diagonalized by using
$\pi_\f(\ga)'\bigwedge\{U_\f(\bp_\mathbb{N})\}'$ as well. Then we get
$$
V_\f=\int^\oplus_{\ce(\cs_{\bp_{\bn}}(\ga))}V_\psi\di\m(\psi)\,.
$$
This implies that the center $\gz_\f$ coincides with its even part  $\gz_\f\bigwedge\{V_\f\}'$.
As in our situation $\gz^\perp_\f=\gz_\f\bigwedge\{V_\f\}'$ (cf. \eqref{veve}), we obtain
\begin{equation}
\label{prinvp}
\gz^\perp_\f=\gz_\f\bigwedge\{V_\f\}'=\gz_\f\,.
\end{equation}
Collecting together Lemma \ref{invper} and \eqref{prinvp}, we get
$$
\gz^\perp_\f=\gz^\perp_\f\bigwedge\{U_\f(\bp_\mathbb{N})\}'=\gz_\f\bigwedge\{U_\f(\bp_\mathbb{N})\}'=\gb_{\bp_\mathbb{N}}(\f)=
\F_\f(\pi_\f(\ga)'')\,.
$$
Let now $X\in\pi_\f(\ga)''$ together with its direct integral decomposition
$X=\int^\oplus_{\ce(\cs_{\bp_{\bn}}(\ga))}X_\psi\di\m(\psi)$. By Kaplansky Density Theorem, there exists a sequence
$\{A_n\}_{n\in\bn}\subset\ga$ such that
$\pi_\f(A_n)\to X$ in the strong operator topology. This implies, by eventually passing to a subsequence, that
$$
\psi(A_n)=\langle\pi_\psi(A_n)\Om_\psi,\Om_\psi\rangle
\to\langle X_\psi\Om_\psi,\Om_\psi\rangle
$$
$\m$--almost everywhere. Fix now $x,y\in\ch_\f$. By Lebesgue Dominated Convergence Theorem, we get
\begin{align*}
\langle\F_\f(X)&x,y\rangle=\lim_n
\langle\F_\f(\pi_\f(A_n))x,y\rangle\\
=&\lim_n\int_{\ce(\cs_{\bp_{\bn}}(\ga))}\langle\pi_\psi(A_n)\Om_\psi,\Om_\psi\rangle
\langle x_\psi,y_\psi\rangle\di\m(\psi)\\
=&\int_{\ce(\cs_{\bp_{\bn}}(\ga))}\langle X_\psi\Om_\psi,\Om_\psi\rangle
\langle x_\psi,y_\psi\rangle\di\m(\psi)\\
=&\left\langle\left(\int^\oplus_{\ce(\cs_{\bp_{\bn}}(\ga))}\langle X_\psi\Om_\psi,\Om_\psi\rangle\idd_{\ch_\psi} \di\m(\psi)\right)x,y\right\rangle\,,
\end{align*}
which leads as a particular case,
\begin{align*}
&\langle X\Om_\f,\Om_\f\rangle=\int_{\ce(\cs_{\bp_{\bn}}(\ga))}\langle X_\psi\Om_\psi,\Om_\psi\rangle\di\m(\psi)\\
=&\int_{\ce(\cs_{\bp_{\bn}}(\ga))}\langle X_\psi\Om_\psi,\Om_\psi\rangle
\langle \Om_\psi,\Om_\psi\rangle\di\m(\psi)=\langle\F_\f(X)\Om_\f,\Om_\f\rangle\,.
\end{align*}
\end{proof}
Once having established in Proposition \ref{havest} that
$\gz^\perp_\f=\gb_{\bp_\mathbb{N}}(\f)$, and $\F_\f$ is the conditional
expectation onto the tail algebra preserving the vector state
$\langle\,{\bf\cdot}\,\Om_\f,\Om_\f\rangle$, it is meaningful to check that any state
$\f\in\cs_{\bp_{\bn}}(\ga)$ is conditionally independent and identically distributed w.r.t. the tail algebra.
\begin{thm}
\label{extdef}
A state $\f\in\cs(\ga)$ is symmetric if and only if the corresponding stochastic process is conditionally independent and identically distributed w.r.t. the tail algebra.
\end{thm}
\begin{proof}
Fix any state $\f\in\ga$. Thanks to $\gz^\perp_\f\subset \gz_\f$ (cf. Theorem 2.6.5 of \cite{BR1}), by Proposition 3.1.1 of \cite{S} we can decompose
$\f=\int_{\cs(\ga)}\psi\di\n(\psi)$ where $\n$ is the orthogonal measure associated to $\gz^\perp_\f$.
By reasoning as in Proposition \ref{havest} (cf. Lemma 8.4.1 of \cite{DX0}), we have $\pi_\f(\ga)''=\int^\oplus_{\cs(\ga)}\pi_\psi(\ga)''\di\n(\psi)$,
and in addition for $X=\int^\oplus_{\cs(\ga)}X_\psi\di\g(\psi)\in\pi_\f(\ga)''$,
the map $E_\f:\pi_\f(\ga)''\to\gz^\perp_\f$ given by
$$
E_\f(X):=\int^\oplus_{\cs(\ga)}\langle X_\psi\Om_\psi,\Om_\psi\rangle \idd_{\ch_\psi}\di\n(\psi)
$$
defines a conditional expectation of $\pi_\f(\ga)''$ onto the tail algebra
$\gz^\perp_\f$ preserving the vector state $\langle\,{\bf\cdot}\,\Om_\f,\Om_\f\rangle$. Suppose that the stochastic process $(\iota_j:\bm_2(\bc)\to\ga)_{j\in \mathbb{N}}$ corresponding to $\f$ is conditionally independent and
identically distributed w.r.t $\gz^\perp_\f$.
By the anticommutation relations, the dense $*$--algebra of the localized elements of $\ga$ coincides with the linear span of terms of the type $\iota_{i_1}(A_1)\iota_{i_2}(A_2)\cdots\iota_{i_n}(A_n)$, where $A_1,\ldots,A_n\in\bm_{2}(\bc)$ and the integers $i_j$ appear only once in the sequence.
By a standard density argument we can reduce the matter to objects of this form. Put 
$A:=\iota_{i_1}(A_1)\iota_{i_2}(A_2)\cdots\iota_{i_n}(A_n)$.
For $g\in\bp_\bn$, we get
$$
\a_g(A)=\iota_{g(i_1)}(A_1))\iota_{g(i_2)}(A_2)\cdots\iota_{g(i_n)}(A_n)\,,
$$
where the $\a_g(i_j)$ again appear only once in the sequence. We compute (cf. Definition \ref{cocaind}),
\begin{align*}
&\f\big(\iota_{i_1}(A_1)\iota_{i_2}(A_2)\cdots\iota_{i_n}(A_n)\big)\\
=&\big\langle E_\f\big(\pi_\f\big(\iota_{i_1}(A_1)\iota_{i_2}(A_2)\cdots\iota_{i_n}(A_n)\big)\big)
\Om_\f,\Om_\f\big\rangle\\
=&\big\langle E_\f\big(\pi_\f\big(\iota_{i_1}(A_1)\big)\big)E_\f\big(\pi_\f\big(\iota_{i_2}(A_2)\big)\big)\cdots
E_\f\big(\pi_\f\big(\iota_{i_n}(A_n)\big)\big)\Om_\f,\Om_\f\big\rangle\\
=&\big\langle E_\f\big(\pi_\f\big(\iota_{g(i_1)}(A_1)\big)\big)E_\f\big(\pi_\f\big(\iota_{g(i_2)}(A_2)\big)\big)\cdots
E_\f\big(\pi_\f\big(\iota_{g(i_n)}(A_n)\big)\big)\Om_\f,\Om_\f\big\rangle\\
=&\big\langle E_\f\big(\pi_\f\big(\iota_{g(i_1)}(A_1)\iota_{g(i_2)}(A_2)\cdots\iota_{g(i_n)}(A_n)\big)
\big)\Om_\f,\Om_\f\big\rangle\\
=&\f\big(\iota_{g(i_1)}(A_1)\iota_{g(i_2)}(A_2)\cdots\iota_{g(i_n)}(A_n)\big)\\
=&\f\circ\a_g\big(\iota_{i_1}(A_1)\iota_{i_2}(A_2)\cdots\iota_{i_n}(A_n)\big)\,,
\end{align*}
that is $\f$ is symmetric.

Let now $\f\in\cs_{\bp_{\bn}}(\ga)$. Then the conditional expectation $\F_\f$ onto the tail algebra $\gz^\perp_\f$
is invariant, i.e. $\F_\f=\F_\f\circ\ad U_\f(g)$, $g\in\bp_\bn$. Thus, the associated process is identically distributed. We now prove that it is conditionally independent w.r.t. $\gz^\perp_\f$. To this aim, for  $I\subset\mathbb{N}$ finite denote by $\gt_{\f,I}$ the  von
Neumann algebra given by $\pi_\f(\carf(I))\bigvee \gz^\perp_\f$.
Fix two finite subsets $I_1$, $I_2$ of $\mathbb{N}$, with
$I_1\cap I_2=\emptyset$, and consider $X_j\in\gt_{\f,I_j}$, $j=1,2$.
We then easily get for such elements,
$$
X_j=\sum_{\b_j }B^{(j)}_{\b_j}\eps^{(j)}_{\b_j}\,,\quad j=1,2\,,
$$
where $\{B^{(j)}_{\b_j }\}\subset\gz^\perp_\f$, and
$\{\eps^{(j)}_{\b_j}\}\subset\pi_\f(\carf(I_j))$ are any system of matrix--units (cf \cite{T1}, Definition IV.1.7) for
$\pi_\f(\carf(I_j))\sim\carf(I_j)$, $j=1,2$. To complete the proof we need to show that
$\Phi_\f(X_1X_2)=\Phi_\f(X_1)\Phi_\f(X_2)$. As $I_1$, $I_2$ are disjoint, and
$\ce(\cs_{\bp_{\bn}}(\ga))$ is made of product states, it easily follows
$\F_\f(R_1R_2)=\F_\f(R_1)\F_\f(R_2)$, whenever $R_j\in\pi_\f(\carf(I_j))$, $j=1,2$. Concerning the general case, we compute,
\begin{align*}
\F_\f(X_1X_2)=&\F_\f\bigg(\sum_{\b_1,\b_2}
B^{(1)}_{\b_1 }
\eps^{(1)}_{\b_1}B^{(2)}_{\b_2}\eps^{(2)}_{\b_2}\bigg)
=\sum_{\b_1,\b_2}
\F_\f\big(B^{(1)}_{\b_1 }
B^{(2)}_{\b_2}\eps^{(1)}_{\b_1}\eps^{(2)}_{\b_2}\big)\\
=&\sum_{\b_1,\b_2}B^{(1)}_{\b_1 }B^{(2)}_{\b_2}
\F_\f\big(\eps^{(1)}_{\b_1}\eps^{(2)}_{\b_2}\big)
=\sum_{\b_1,\b_2}B^{(1)}_{\b_1}B^{(2)}_{\b_2}
\F_\f\big(\eps^{(1)}_{\b_1}\big)\F_\f\big(\eps^{(2)}_{\b_2}\big)\\
=&\bigg(\sum_{\b_1}B^{(1)}_{\b_1 }\F_\f\big(\eps^{(1)}_{\b_1}\big)\bigg)
\bigg(\sum_{\b_2}B^{(2)}_{\a_2,\b_2 }\F_\f\big(\eps^{(2)}_{\b_2}\big)\bigg)\\
=&\bigg(\sum_{\b_1}\F_\f\big(B^{(1)}_{\b_1 }\eps^{(1)}_{\b_1}\big)\bigg)
\bigg(\sum_{\b_2}\F_\f\big(B^{(2)}_{\b_2 }\eps^{(2)}_{\b_2}\big)\bigg)\\
=&\F_\f(X_1)\F_\f(X_2)\,.
\end{align*}
\end{proof}
\begin{rem}
\label{inften}
By following the same lines of the above proof, we can show that a stochastic process factoring through the countably infinite tensor product of a single separable $C^*$--algebra is conditionally independent and identically distributed w.r.t. the tail algebra if and only if the corresponding state is symmetric.
\end{rem}

\section{symmetric states in free and $q$--deformed probability}
\label{sec4}

The $q$--deformed Commutation Relations were introduced in Quantum Physics in \cite{FB}. The reader is referred to \cite{DyF} and the references cited therein for further details.
Suppose $-1<q<1$ and fix an Hilbert space $\ch$. The $q$--deformed Fock space $\G_q(\ch)$ is the completion of
the algebraic linear span of the vacuum vector $\Om$, together with
vectors
$$
f_1\otimes\cdots\otimes f_n\,,\quad
f_j\in\ch\,,j=1,\dots,n\,,n=1,2,\dots
$$
with respect to the $q$--deformed inner product
\begin{equation*}
\langle f_1\otimes\cdots\otimes
f_n\,,g_1\otimes\cdots\otimes g_m\rangle_q
:=\d_{n,m}\sum_{\pi\in\bp_n}q^{i(\pi)}\langle
f_1\,,g_{\pi(1)}\rangle_\ch\cdots\langle f_n\,,g_{\pi(n)}\rangle_\ch\,,
\end{equation*}
$\d_{n,m}$ being the Kronecker symbol, $\bp_n$ the symmetric group of $n$ elements, and $i(\pi)$ the
number of inversions of $\pi\in\bp_n$.
Fix $f,f_1,\ldots,f_n\in\ch$. Define the creator $a_q^\dagger(f)$ as
$$
a_q^\dagger(f)\Omega=f\,,\quad a_q^\dagger(f)f_1\otimes\cdots\otimes f_n=f\otimes
f_1\otimes\cdots\otimes f_n\,,
$$
and the annihilator $a(f)$ as
\begin{align*}
%label{qann}
&a_q(f)\Omega=0\,,\quad
a_q(f)(f_1\otimes\cdots\otimes f_n)\\
=&\sum_{k=1}^nq^{k-1}\langle
f_k,f\rangle_\ch f_1\otimes\cdots f_{k-1}\otimes
f_{k+1}\otimes\cdots\otimes f_n\,.
\end{align*}
$a_q^\dagger(f)$ and $a_q(f)$ are mutually adjoint with respect to the $q$--deformed inner
product, and satisfy the commutation relations
$$
a_q(f)a_q^\dagger(g)-qa_q^\dagger(g)a_q(f)=\langle
g,f\rangle_{\ch}\idd\,,\qquad f,g\in\ch\,.
$$
The limiting cases are the Canonical Commutation Relations (Bosons) when $q=1$, and the
Canonical Anticommutation Relations (Fermions) for $q=-1$, the latter treated exhaustively in \cite{CrF} and in Section \ref{sec3}. The case $q=0$ corresponds to the free group reduced $C^*$--algebra (see below).

The concrete $C^*$--algebra $\gar_q$ and its subalgebra $\gg_q$, acting on $\G_q(\ch)$ are the unital
$C^*$--algebras generated by the annihilators $\{a(f)\mid f\in\ch\}$, and the selfadjoint part of
annihilators $\{s_q(f)\mid f\in\ch\}$,
$$
s_q(f):=a_q(f)+a_q^\dagger(f)\,,\quad f\in\ch\,,
$$
respectively. The Fock vacuum expectation is the state on both the mentioned $C^*$--algebras defined as
$\om_q:=\langle\,{\bf\cdot}\,\Om,\Om\rangle$.
Notice that $\gar_0$ is an extension of
the Cuntz algebra $\co_n$ by the compact operators if $\ch$ has dimension $2\leq n<\infty$, see e.g. see \cite{VDN}, pag. 6.
\begin{rem}
\label{cuntz}
The Cuntz algebra $\co_\infty$ coincides with $\gar_0$ with $\ch$ separable infinite dimensional. Its generators are
$\big\{a_0^\dagger(e_j)\mid j\in\bn\big\}$, where $\{e_j\}_{j\in\bn}$ is any orthonormal basis of $\ch$.
\end{rem}

In order to study the symmetric states on $\gar_q$ and $\gg_q$,
we need to consider the action of the group of permutations $\bp_J$ on them. As in the previous sections, we get
$J\sim \mathbb{Z}$. Then $\gar_q$ and $\gg_q$ are concrete $C^*$--algebras on $\G_q(\ell^2(\bz))$.
If $i\in\mathbb{Z}$ and $e_i\in\ell^2(\bz)$ is the sequence taking value $1$ at $i$ and zero elsewhere,
we denote $a_i:=a(e_i)$, $a_i^\dagger:=a^\dagger(e_i)$ and $s_i:=a_i+a^\dagger_i$. As usual, the group of permutations
$\bp_\bz$ naturally acts on $\gar_q$ and $\gg_q$ by
$$
\a_g(a_i):=a_{g(i)}\,,\quad i\in\bz\,,\,g\in\bp_\bz\,,
$$
and the Fock vacuum is invariant under such an action. The group $\bz$ also acts on both
$\gar_q$ and $\gg_q$ as powers of the right shift $\b$, uniquely determined by $\b(a_i):=a_{i+1}$, $i\in\bz$,
on the generators. In addition, the Fock vacuum $\om_q$ is shift--invariant, and it is shown in \cite{DyF} that it is the unique invariant state. The same state is moreover the unique symmetric one on both $\gar_q$ and $\gg_q$, as the following proposition  shows.
\begin{prop}
\label{frepp}
The set $\cs_{\bp_\bz}(\gar_q)$ and $\cs_{\bp_\bz}(\gg_q)$ of the invariant states under the action of $\bp_\bz$ on $\gar_q$ and
$\gg_q$ respectively, consist of a singleton which is precisely the Fock vacuum expectation.
\end{prop}
\begin{proof}
Fix the integers $k<l$ and take an arbitrary polynomial
$A$ in $a_q(f)$ and $a_q^\dagger(f)$, or in $s_q(f)$ where the degree zero corresponds to the multiple of the identity, with
$f\in\ell^2([k,l])\subset \ell^2(\bz)$. By Lemma \ref{permtra}, there is a cycle $g_A\in\bp_\bz$ such that
$\b(A)=\a_{g_A}(A)$. Thus,
$$
\f(\b(A))=\f(\a_{g_A}(A))=\f(A)\,,
$$
provided that
$\f\in\cs_{\bp_\bz}(\gar_q)$ or $\f\in\cs_{\bp_\bz}(\gg_q)$, respectively. As the polynomials as above are dense in
$\gar_q$ or $\gg_q$ respectively,
we conclude that, if $\f$ is symmetric, then it is shift--invariant. But it is shown in Theorem 3.3 and Corollary 3.4 of \cite{DyF}, that the Fock vacuum expectation is the unique shift invariant state.
\end{proof}
An immediate application of Proposition \ref{frepp} gives that the Fock vacuum $\om_0$ is the unique symmetric state on the Cuntz algebra $\co_\infty$ (see Remark \ref{cuntz}), since the natural action of the permutation group $\bp_\infty$ on this algebra.

Let now $\bbf_\infty$ be the free group with countably many generators $\{g_i\mid i\in\mathbb{Z}\}$. For
$i_1\neq i_2\neq\cdots\neq i_n$ where the indices can appear more than once in the string, and
$k_1,k_2,\dots,k_n\in\bz\backslash\{0\}$, denote $w:=g_{i_1}^{k_1}\cdots g_{i_n}^{k_n}$ a (reduced) word.
The length of $w\in\bbf_\infty$ is defined as
$$
w:=|k_1|+|k_2|+\cdots+|k_n|\,.
$$
The empty word $w_{\emptyset}$, whose length is zero by definition, is the unity of $\bbf_\infty$. All the words of arbitrary length, together with the natural group operations, generate $\bbf_\infty$. The universal $C^*$--algebra generated by the linear combinations of
$\{\l_w\mid w\in\bbf_\infty\}$, together with the relations $\l_w^*=\l_{w^{-1}}$, $w\in\bbf_\infty$, is the group $C^*$--algebra
$C^*(\bbf_\infty)$. This is nothing but the unital free product $C^*$--algebra of the group $\bz$. The concrete $C^*$--algebra $C_r^*(\bbf_\infty)$ generated by (left) regular representation is precisely the
reduced group $C^*$--algebra $C_r^*(\bbf_\infty)$. It is well--known that it differs from $C^*(\bbf_\infty)$ as $\bbf_\infty$ is not amenable. The GNS representation $\pi_\t$ associated to the tracial state $\t$ on $C^*(\bbf_\infty)$, uniquely defined as
$$
\t\bigg(\sum_{w\in\bbf_\infty}a_w\l_w\bigg):=a_{w_{\emptyset}}\,,
$$
generates the left regular representation, so $\pi_\t\big(C^*(\bbf_\infty)\big)=C_r^*(\bbf_\infty)$. Finally the group $\bp_\infty$ acts in a natural way on $C^*(\bbf_\infty)$ and
$C_r^*(\bbf_\infty)$.
\begin{cor}
\label{cofr}
The tracial state $\t$ is the unique symmetric state on $C_r^*(\bbf_\infty)$.
\end{cor}
\begin{proof}
The proof directly follows collecting together Theorem 2.6.2 of \cite{VDN}, and Proposition \ref{frepp}.
\end{proof}
The last part of the section is devoted to determine some remarkable ergodic properties of the Haagerup states.
Recall that the Haagerup states on $C^*(\bbf_\infty)$, labelled by $\l\in(0,+\infty)$, are defined in \cite{Ha} as
$$
\f_\l(w):=e^{-\l|w|}\,.
$$
The case $\l=+\infty$ corresponds to the tracial state and it is covered by Corollary \ref{cofr}.
It generates the regular representation of $\bbf_\infty$ and in addition, it is shown in Corollary 3.2 of \cite{Ha} that is the unique Haagerup state
which is normal w.r.t. the regular representation. Namely, it is the unique state of the family whose corresponding stochastic process
take values in the reduced free product $C^*$--algebra $C^*_r(\bbf_\infty)$. The Haagerup states are automatically symmetric by construction, and satisfy the product state condition
$\f_\l(vw)=\f_\l(v)\f_\l(w)$ if $I_v\cap I_w=\emptyset$.\footnote{The Haagerup states satisfy a stronger
condition $\f_\l(vw)=\f_\l(v)\f_\l(w)$ if $|vw|=|v|+|w|$.} But they do not fulfil the block--singleton condition for $\l\in(0,+\infty)$. In fact, if $i\neq j$,
$$
\f_\l(g_ig_jg_i^{-1})=e^{-3\l}\neq e^{-\l}=\f_\l(g_j)\f_\l(\idd)=\f_\l(g_j)\f_\l(g_ig_i^{-1})\,.
$$
By (i) of Theorem \ref{frecazo} they are weakly clustering, then extremal (i.e. ergodic under the action of $\bp_\infty$)
thanks to Proposition 3.1.10. But (ii) of Theorem \ref{frecazo} gives that
$\pi_{\f_\l}\big(C^*(\bbf_\infty)\big)''\bigwedge \big\{U_{\f_\l}(\bp_\infty)\big\}'$ cannot be trivial, this property being crucial for the proof of the following
\begin{thm}
%\label{center}
For $\l\in(0,+\infty)$, the support $s(\f_\l)\in C^*(\bbf_\infty)^{**}$ of $\f_\l$ does not belong to $Z(C^*(\bbf_\infty)^{**})$.
\end{thm}
\begin{proof}
As we have shown, the Haagerup states are extremal symmetric. Suppose that for the support of $\f_\l$ in the bidual is central, then the cyclic vector $\Om_{\f_\l}$ is also separating for $\pi_{\f_\l}\big(C^*(\bbf_\infty)\big)''$. By Theorem 4.3.20 in \cite{BR1},
it follows that $\pi_{\f_\l}\big(C^*(\bbf_\infty)\big)''\bigwedge \big\{U_{\f_\l}(\bp_\infty)\big\}'=\bc I$. As a consequence, by Theorem \ref{frecazo}, $\f_\l$  satisfies
the block--singleton condition, This contradicts the above discussion.
\end{proof}
In conclusion, as realized throughout the section, the structure of the convex set made by the symmetric states and its boundary, change radically when one considers
$C^*_r(\bbf_\infty)$ or $C^*(\bbf_\infty)$, passing from a singleton to a richer structure which contains, perhaps properly, all the Haagerup states.

\section{the boolean case}
\label{sec5}
Let $\ch$ be a complex Hilbert space. Recall that the Boolean Fock space over $\ch$ (cf. \cite{BGS}) is given by $\G(\ch):=\mathbb{C}\oplus \ch$,
where the vacuum vector $\Om$ is $(1,0)$.
On $\Gamma(\ch)$ we define the creation and annihilation operators, respectively given for $f\in \ch$, by
$$
b^\dagger(f)(\alpha\oplus g):=0\oplus \alpha f,\,\,\,\, b(f)(\alpha\oplus g):=\langle g,f\rangle_\ch \oplus 0,\,\,\, \alpha\in\mathbb{C},\, g\in\ch.
$$
They are mutually adjoint, and satisfy the following relations for $f,g\in \ch$,
$$
b(f)b^\dagger(g)=\langle g, f\rangle_\ch \langle\,{\bf\cdot}\,,\Om\rangle\Om\,,\quad
b^\dagger(f)b(g)=\langle\,{\bf\cdot}\,,0\oplus g\rangle 0\oplus f\,.
$$
As in Section \ref{sec4}, we consider the unital $C^*$--algebras acting on $\Gamma(\ch)$, which are respectively
generated by the annihilators $\{b(f)\mid f\in\ch\}$, and the selfadjoint part of annihilators $\{r(f)\mid f\in\ch\}$, where
$$
r(f):=b(f)+b^\dagger(f)\,,\quad f\in\ch\,.
$$
Again, we are interested in the action of the group of permutations $\bp_J$ on these algebras. Then, as usual, we put $J\sim \mathbb{Z}$, and take $\ch=\ell^2(\mathbb{Z})$. From now on, we denote the vacuum vector by $e_\#$, and it can
be seen as an element of $\ell^2(\{\#\}\cup\mathbb{Z})=\Gamma(l^2(\mathbb{Z}))$. As a consequence, $\{e_\#, e_i|i\in \mathbb{Z}\}$ is an orthonormal basis for this space,
where for each $i\in \mathbb{Z}$, $e_i$ is the sequence taking value $1$ on $i$ and $0$ elsewhere.
For each $j\in \mathbb{Z}$, define $b_j:=b(e_j)$, $b^\dagger_j:=b^\dagger(e_j)$, and $r_j:=r(e_j)$, $j\in\bz$. With
$\cf(\ell^2(\{\#\}\cup\mathbb{Z}))$ and $\ck(\ell^2(\{\#\}\cup\mathbb{Z}))$ we will denote respectively,
the finite rank operators and the compact linear operators acting on $\ell^2(\{\#\}\cup\mathbb{Z})$.
\begin{prop}
\label{bolca}
The unital $C^*$--algebras acting on the Boolean Fock space $\Gamma(\ell^2(\bz))$ generated by the Boolean annihilators
$\{b_j\mid j\in\bz\}$, or by their selfadjoint parts
$\{r_j\mid j\in\bz\}$,
are equal and coincide with $\ck(\ell^2(\{\#\}\cup\bz))+\bc\idd$.
\end{prop}
\begin{proof}
First we note that
$\alg\{b_j\mid j\in\bz\}$ generates all of $\cf(\ell^2(\{\#\}\cup\bz))$. In fact, consider in $\cb(\ell^2(\{\#\}\cup\bz))$, the canonical system of matrix--units $\{\varepsilon_{ij}\mid i,j\in(\{\#\}\cup\bz)\}$. It is easy to check that
$$
\varepsilon_{\#j}=b_j\,,\,\, \varepsilon_{j\#}=b^\dagger_j\,,\,\,
\varepsilon_{\#\#}=b_ib^\dagger_i\,,\quad \varepsilon_{ij}=b^\dagger_ib_j\,,\quad i,j\in\bz\,,
$$
that is the assertion. The equality
$\alg\{r_j\mid j\in\bz\}=\cf(\ell^2(\{\#\}\cup\bz))$ follows after noticing that
$$
r_ir_j=\d_{ij}\varepsilon_{\#\#}+\varepsilon_{ij}\,,
r^2_i-r^2_ir^2_j=\varepsilon_{ii}-\d_{ij}\varepsilon_{ij}\quad i,j\in\bz\,.
$$
Since $\ck(\ell^2(\{\#\}\cup\bz))$ is the norm closure of $\cf(\ell^2(\{\#\}\cup\bz))$ in $\cb(\ell^2(\{\#\}\cup\bz))$, the thesis follows.
\end{proof}
It is worth noticing that also the universal $C^*$--algebra generated by Boolean annihilators, or equivalently by their selfadjoint part, is isomorphic
to the algebra of compact operator on $\ell^2(\{\#\}\cup\bz)$, see e.g. the arguments outlined in \cite{Sz}.

Denote $\gpb=\ck(\ell^2(\{\#\}\cup\bz))+\bc I$ the unital Boolean $C^*$--algebra. The finite permutations and the shift naturally
act on the indices of $\bz$, leaving invariant the vacuum index $\#$ corresponding to the vacuum vector $e_\#$. Thus, we have
natural actions of the permutations and the shift on $\gpb$ and
$\ck(\ell^2(\{\#\}\cup\bz))$, both denoted by an abuse of notations, by
$\{\a_g\mid g\in\bp_\bz\}$ and $\b$, respectively.
\begin{lem}
\label{shiv}
The unique invariant state under the shift for $\ck(\ell^2(\{\#\}\cup\bz))$ is the vacuum state.
\end{lem}
\begin{proof}
Let $U$ be the canonical implementation of the shift on the Boolean Fock space:
$$
Ue_\#=e_\#\,,\quad Ue_k=e_{k+1}\,,\,\,\,k\in\bz\,.
$$
Fix any state $\om_T\in\cs(\ck(\ell^2(\{\#\}\cup\bz)))$, together with its representation through a positive trace class operator
$$
T=\sum_{\l\in\s_{{\rm pp}}(T)}\l E_\l\,
$$
where "$\s_{{\rm pp}}$" stands for pure points spectrum.
The fact that $\om_T$ is invariant w.r.t. the shift leads to $T=UTU^*$, which turns out to be equivalent to
$$
U=\sum_{\l\in\s_{{\rm pp}}(T)}E_\l UE_\l\,.
$$
As the unique eigenspace of $U$ is $\bc e_\#$ (see e.g. Section 6 of \cite{F}), the last condition is possible only if
$T=\langle\,{\bf\cdot}\, e_\#,e_\#\rangle e_\#$.
\end{proof}
\begin{prop}
We have for the compact convex set of the symmetric states,
$$
\cs_{\bp_{\bz}}(\gpb)=\{\g\om_\#+(1-\g)\om_\infty\mid\g\in[0,1]\}\,,
$$
where $\om_\#=\langle\,{\bf\cdot}\,e_\#,e_\#\rangle$ is the Fock vacuum state, and
$$
\om_\infty(A+aI):=a\,,\quad A\in \ck(\ell^2(\{\#\}\cup\bz))\,, a\in\bc\,.
$$
\end{prop}
\begin{proof}
Fix $\om\in\cs_{\bp_{\bz}}(\gpb)$. Arguing as in Proposition \ref{frepp}, one finds that $\om$ is also shift--invariant. Its
restriction $\om\lceil_{\ck(\ell^2(\{\#\}\cup\bz))}$ yields a positive functional which is also shift--invariant. Thus, by Proposition \ref{shiv} we get,
$$
\om\lceil_{\ck(\ell^2(\{\#\}\cup\bz))}
=\big\|\om\lceil_{\ck(\ell^2(\{\#\}\cup\bz))}\big\|\om_\#\,.
$$
This means that $\om=\g\om_\#+(1-\g)\om_\infty$, where $\g=\big\|\om\lceil_{\ck(\ell^2(\{\#\}\cup\bz))}\big\|$.
\end{proof}
To end the present section, we show that there are plenty of Boolean processes for which the tail algebra is not expected, that is no conditional expectation onto such an algebra which preserves the state corresponding to the process under consideration.
To simplify, we consider a pure state $\om_\xi:=\langle\,{\bf\cdot}\,\xi,\xi\rangle$, $\xi\in\ell^2(\{\#\}\cup\bz)$ being a unit vector. In this situation $\pi_{\om_\xi}(\gpb)''=\cb(\ell^2(\{\#\}\cup\bz))$. Concerning the tail algebra, we get with $I_n:=\{k\in\bz\mid |k|>n\}$ and 
$P_\#:=\langle\,{\bf\cdot}\, e_\#,e_\#\rangle e_\#$,
$$
\gz^\perp_{\om_\xi}=\bigcap_{n\in\bn}\cb(\ell^2(\{\#\}\cup I_n))\bigoplus\bc P_{\ell^2(\bz\backslash I_n)}
=\bc P_\#\oplus\bc P_\#^\perp\,,
$$
where $\cb(\ell^2(\{\#\}\cup I_n))$ is considered as a non unital subalgebra of 
$\cb(\ell^2(\{\#\}\cup\bz))$ in a canonical way.
Notice that each conditional expectation $F$ onto $\gz^\perp_{\om_\xi}$ satisfies $F=F\circ E$, with
$$
E(A)=\om_\#(A)P_\#+P_\#^\perp AP_\#^\perp\,,\quad A\in \cb(\ell^2(\{\#\}\cup \bz))\,.
$$
By taking into account that any conditional expectation of a $C^*$--algebra with unity $\idd$, into $\bc\idd$ is just given by a state, we conclude that any conditional expectation $F=F_\f$ as above 
assumes the form
$$
F_\f(A)=\om_\#(A)P_\#+\f(P_\#^\perp AP_\#^\perp)P_\#^\perp\,,\quad A\in \cb(\ell^2(\{\#\}\cup \bz))\,,
$$ 
where $\f$ is any state, not necessarily normal, on $\cb(\ell^2(\bz))$, the last viewed again as a non unital subalgebra of 
$\cb(\ell^2(\{\#\}\cup\bz))$.
We now show that $\gz^\perp_{\om_\xi}$ cannot be expected if $0<|\langle e_\#,\xi\rangle|<1$. For the rank--one operator
$A=\langle \,{\bf\cdot}\,,\xi\rangle e_\#$ and any $\f$ as above, we get 
$$
\frac{\om_\xi(F_\f(A))}{\om_\xi(A)}=|\langle e_\#,\xi\rangle|^2\,.
$$
The examples of Boolean stochastic processes described above explain that, contrarily to the classical case,
the condition to be independent and identically distributed w.r.t. the tail algebra (cf. Definition \ref{cocaind}), cannot be formulated in the general case,
without mentioning the {\it a--priori} existence of a preserving conditional expectation.

\end{document}